# ORDER SELECTION FOR SAME-REALIZATION PREDICTIONS IN AUTOREGRESSIVE PROCESSES

By Ching-Kang Ing and Ching-Zong Wei[1]

*Academia Sinica and National Taiwan University*

Assume that observations are generated from an infinite-order autoregressive [AR($\infty$)] process. Shibata [*Ann. Statist.* **8** (1980) 147–164] considered the problem of choosing a finite-order AR model, allowing the order to become infinite as the number of observations does in order to obtain a better approximation. He showed that, for the purpose of predicting the future of an independent replicate, Akaike's information criterion (AIC) and its variants are asymptotically efficient. Although Shibata's concept of asymptotic efficiency has been widely accepted in the literature, it is not a natural property for time series analysis. This is because when new observations of a time series become available, they are not independent of the previous data. To overcome this difficulty, in this paper we focus on order selection for forecasting the future of an observed time series, referred to as same-realization prediction. We present the first theoretical verification that AIC and its variants are still asymptotically efficient (in the sense defined in Section 4) for same-realization predictions. To obtain this result, a technical condition, easily met in common practice, is introduced to simplify the complicated dependent structures among the selected orders, estimated parameters and future observations. In addition, a simulation study is conducted to illustrate the practical implications of AIC. This study shows that AIC also yields a satisfactory same-realization prediction in finite samples. On the other hand, a limitation of AIC in same-realization settings is pointed out. It is interesting to note that this limitation of AIC does not exist for corresponding independent cases.

**1. Introduction.** To select a model for the realization of a stationary time series, it is common to assume that the realization comes from an au-toregressive moving-average (ARMA) process whose AR and MA orders are known to lie within prescribed finite intervals. Then a model selection procedure is used to select orders within these intervals and thereby determine a model for the data. However, as pointed out by Burnham and Anderson [9], it is not common for the true model to be a function of a small number of unknown parameters, and a model having many parameters is sometimes essential to obtain a better approximation of the true model. From this perspective, a more flexible alternative to the ARMA assumption is the assumption that data are generated by an AR($\infty$) process. In this situation, the focus of model selection is usually placed on the forecasting ability of the chosen model, and not on the correctness of the selection.

Shibata [27] gave the first justification for several model selection criteria along this line. He considered the problem of choosing a finite-order AR model, allowing the order to become infinite as the number of observations does. He showed that for the purpose of forecasting the future of an independent replicate, which is referred to as *independent-realization prediction* [see (1.4)], Akaike's information criterion (AIC) [2], the final prediction error (FPE) method [1] and $S_n(k)$ [27] are asymptotically efficient in the sense that no other selection criterion achieves a smaller limiting mean square prediction error criterion value. (Since this is an asymptotic result, the name AIC could also be thought of as an acronym for "Asymptotic Information Criterion.") Based on a similar analysis, Bhansali [5] extended Shibata's result to the case of multistep predictions. However, Shibata's concept of asymptotic efficiency, which focuses on independent-realization predictions, is not a natural property for time series analysis, because when new observations of a time series become available, they are usually dependent on the previous data. So far, no time series model selection theory has been established without this unnatural assumption. This motivated our study.

To begin with, let us assume that observations $x_1, \ldots, x_n$ come from a stationary AR($\infty$) process $\{x_t\}$ with

$$(1.1) \qquad x_t + \sum_{i=1}^{\infty} a_i x_{t-i} = e_t, \qquad t = \ldots, -1, 0, 1, \ldots,$$

where $e_t$ is a sequence of independent random noise values with zero mean and variance $\sigma^2$, and the coefficients $a_i$ are absolutely summable. For predicting $x_{n+h}$, $h \geq 1$, we consider the finite-order approximation models AR(1), ..., AR($K_n$). Here, we allow the maximal order, $K_n$, to increase to infinity with $n$ in order to reduce approximation errors. The prediction for $x_{n+h}$ is referred to as *same-realization* prediction. For brevity, our theoretical discussion only focuses on the one-step prediction case, $h = 1$. But the related extensions to cases $h > 1$ are straightforward as discussed in Section 6. When model

AR($k$), $1 \leq k \leq K_n$, is adopted, we use $\hat{\mathbf{a}}_n(k)$ to estimate the model's coefficient vector and use

$$\hat{x}_{n+1}(k) = -\mathbf{x}'_n(k)\hat{\mathbf{a}}_n(k) \tag{1.2}$$

to predict $x_{n+1}$, where $\mathbf{x}_j(k) = (x_j, \ldots, x_{j-k+1})'$ and

$$\hat{\mathbf{a}}_n(k) = (\hat{a}_{1,n}(k), \ldots, \hat{a}_{k,n}(k))'$$

satisfies

$$-\hat{R}_n(k)\hat{\mathbf{a}}_n(k) = \frac{1}{N} \sum_{j=K_n}^{n-1} \mathbf{x}_j(k) x_{j+1}, \tag{1.3}$$

with $N = n - K_n$ and

$$\hat{R}_n(k) = \frac{1}{N} \sum_{j=K_n}^{n-1} \mathbf{x}_j(k) \mathbf{x}'_j(k).$$

Since the difference between $\hat{\mathbf{a}}_n(k)$ and the least squares estimator $\hat{\mathbf{a}}^L_n(k)$, where

$$\hat{\mathbf{a}}^L_n(k) = -\left(\sum_{j=k}^{n-1} \mathbf{x}_j(k)\mathbf{x}'_j(k)\right)^{-1} \sum_{j=k}^{n-1} \mathbf{x}_j(k) x_{j+1},$$

is asymptotically negligible under the assumptions on $K_n$ and $x_t$ we use herein (see Section 2), $\hat{x}_{n+1}(k)$ is still called the least squares predictor. For assessing the model's predictive ability, we consider the second-order (unconditional) mean-squared prediction error (MSPE), $l_n(k)$, of $\hat{x}_{n+1}(k)$, where

$$l_n(k) = E(x_{n+1} - \hat{x}_{n+1}(k))^2 - \sigma^2.$$

In Section 2 some asymptotic properties of $l_n(k)$ from a companion paper [17] are introduced. In particular, Proposition 2 of Section 2 shows that $l_n(k)$ can be uniformly (in $k$) approximated by $L_n(k) = (k/N)\sigma^2 + \|\mathbf{a} - \mathbf{a}(k)\|_R^2$, where $\mathbf{a}(k)$ is defined after (2.3), $\mathbf{a} = (a_1, a_2, \ldots)'$ is an infinite-dimensional vector with $a_i$'s defined in (1.1), and $\|\mathbf{a} - \mathbf{a}(k)\|_R^2$ is defined after (2.6). The first term of $L_n(k)$, $(k/N)\sigma^2$, which is proportional to the order of the candidate model, $k$, can be viewed as a measure of model complexity. The second term of $L_n(k)$, $\|\mathbf{a} - \mathbf{a}(k)\|_R^2$, which decreases as $k$ increases, measures the goodness of fit. Proposition 3 of Section 2 further points out that $L_n(k)$ can also be used to uniformly approximate $l_{n,0}(k) = E(y_{n+1} - \hat{y}_{n+1}(k))^2 - \sigma^2$, the second-order unconditional MSPE for independent-realization predictions. Here $\{y_1, \ldots, y_n\}$ is a realization from an independent copy of $\{x_t\}$, $y_{n+1}$ is the future observation to be predicted, and the predictor $\hat{y}_{n+1}(k)$ is given by

$$\hat{y}_{n+1}(k) = -\mathbf{y}'_n(k)\hat{\mathbf{a}}_n(k), \tag{1.4}$$

with $\hat{\mathbf{a}}_n(k)$ [see (1.3)] obtained from $x_1, \ldots, x_n$ and $\mathbf{y}'_n(k) = (y_n, \ldots, y_{n+1-k})$. Therefore, these two types of MSPEs are asymptotically equivalent. To us this equivalence is somewhat surprising because some recent studies have shown that this equivalence does not hold in other situations; see the discussion after Proposition 3 for details. It can be erroneous to directly assume that the results from same-realization predictions will be the same as those for corresponding independent cases without theoretical justification.

When the order of the least squares predictor is selected by an order selection criterion, due to more complicated probabilistic structures, analyzing the predictor's MSPE becomes more difficult. Section 3 is devoted to this problem. For independent-realization predictions, Theorem 1 of Section 3 provides an asymptotic expression for the second-order unconditional MSPE of $\hat{y}_{n+1}(\hat{k}_n)$,

$$l_{n,0}(\hat{k}_n) = E(y_{n+1} - \hat{y}_{n+1}(\hat{k}_n))^2 - \sigma^2, \tag{1.5}$$

where $1 \leq \hat{k}_n = \hat{k}_n(x_1, \ldots, x_n) \leq K_n$ is an order determined by AIC, FPE, $C_p$ [22], $S_p$ [14] or $S_n(k)$. The reason why $C_p$ and $S_p$ are included in our analysis is given in Remark 2 of Section 3. We are interested in the other criteria because their asymptotic optimalities for independent-realization predictions were justified by Shibata [27] through a conditional version of $l_{n,0}(\hat{k}_n)$, namely,

$$E\{(y_{n+1} - \hat{y}_{n+1}(\hat{k}_n))^2 | x_1, \ldots, x_n\} - \sigma^2;$$

see (3.1) and (4.11) for more details. However, since this paper focuses on the unconditional MSPE, an extension of Shibata's result to the unconditional case is needed. It should be noted that this extension is nontrivial since there are several technical gaps to be bridged, as detailed in Section 5. According to Theorem 1, $\ell_{n,0}(\hat{k}_n)$ with $\hat{k}_n$ selected by these criteria can ultimately achieve the best compromise between model complexity and goodness of fit, provided $\{x_t\}$ is truly an infinite-order AR process. Viewing this result, it is interesting to ask whether AIC [$C_p$, $S_p$, FPE or $S_n(k)$] still possesses a similar property for same-realization predictions. The main difficulty of this question lies in the fact that the *selected orders*, *estimated parameters* and *future observations* are all stochastically dependent in the same-realization case. Since, as observed in (1.5), the future observations are independent of the estimated parameters and the selected order for independent-realization predictions, the approaches used in [27] and Theorem 1 are no longer applicable. To overcome this difficulty, an assumption for $L_n(k)$, assumption (K.6), is introduced in this section. Two examples are given to illustrate that assumption (K.6) is easily met in common practice. Based on this assumption (among others), Theorem 2 (also in Section 3) shows that

$$E(x_{n+1} - \hat{x}_{n+1}(\hat{k}_n))^2 - \sigma^2$$

and

$$E(y_{n+1} - \hat{y}_{n+1}(\hat{k}_n))^2 - \sigma^2,$$

with $\hat{k}_n$ selected by AIC [$C_p$, $S_p$, FPE or $S_n(k)$], have the same asymptotic expressions. Moreover, we also apply the same techniques to analyze some other AIC-like criteria having different penalty functions; see Corollary 1 and Remark 3 of Section 3. Armed with Corollary 1, the performances of these criteria are first evaluated from the same-realization prediction point of view.

In Section 4 the results obtained in Section 3 are re-examined in greater depth. In particular, we show that, for same-realization predictions, the predictor with an order determined by AIC [FPE, $S_p$, $C_p$ or $S_n(k)$] is ultimately no worse than the best predictor among the candidate predictors, $\{\hat{x}_{n+1}(1), \ldots, \hat{x}_{n+1}(K_n)\}$. This property is referred to as asymptotic efficiency; see (4.1). To the authors' knowledge, this is the first result that confirms AIC's (and its variants') validity in same-realization settings. In addition, a simulation study is conducted to illustrate the practical implications of AIC. This study shows that AIC also yields a satisfactory same-realization prediction in many finite-sample situations; see Table 1 in Section 4 for more details. On the other hand, a limitation of AIC in same-realization settings is demonstrated. Empirical results, given in Table 2 in Section 4, reveal that it seems very difficult for AIC to possess strong asymptotic efficiency; see (4.5) for the definition. This is a somewhat interesting discovery because we show at the end of Section 4 that AIC has no such difficulty when it is used for independent-realization predictions. It is worth noting that AIC's asymptotic efficiency is established under the assumption that the underlying process is truly an AR($\infty$) process. If the order of the true model is finite, then the BIC-like criterion, for example, BIC [24] and HQ [13], can choose the smallest true model with probability tending to 1, but AIC does not possess this optimal property (see [26]). Therefore, to achieve optimal same-realization predictions in situations where the underlying AR model has a possibly finite order, further investigation is still required. For ease of reading, all proofs of the results in Section 3 are deferred to Section 5. Concluding remarks are given in Section 6. Discussions of moment restrictions, connections between time series and regression model selections, and extensions to the multivariate case are also given in this section.

**2. Preliminary results.** We first list several assumptions essential to the following analysis.

(K.1) Let $\{x_t\}$ be a linear process satisfying (1.1) with $A(z) = 1 + a_1 z + a_2 z^2 + \cdots \neq 0$ for $|z| \leq 1$. Furthermore, let the coefficients $\{a_i\}$ obey one of the following restrictions: (a) $\sum_{i=1}^{\infty} |a_i| < \infty$, (b) $\sum_{i=1}^{\infty} |i^{1/2} a_i| < \infty$, or (c) $\sum_{i=1}^{\infty} |i a_i| < \infty$.

(K.2) Let the distribution function of $e_t$ be denoted by $F_t$. Some positive numbers $\alpha, \delta$ and $C_0$ exist such that, for all $t = \ldots, -1, 0, 1, \ldots$ and $|x - y| < \delta$,
$$|F_t(x) - F_t(y)| \leq C_0 |x - y|^{\alpha}.$$

(K.3) $\sup_{-\infty < t < \infty} E|e_t|^s < \infty$, $s = 1, 2, \ldots$.

(K.4) Let $K_n$ be chosen to satisfy
$$C_l \leq \frac{K_n^{2+\delta_1}}{n} \leq C_u$$
for some positive numbers $\delta_1$, $C_l$ and $C_u$.

(K.5) $a_n \neq 0$ for infinitely many $n$.

REMARK 1. (K.1)(a) implies that $x_t$ has a one-sided infinite moving-average representation ([30], page 245),

$$x_t = \sum_{i=0}^{\infty} b_i e_{t-i}, \tag{2.1}$$

where the $b_i$ are absolutely summable with $b_0 = 1$, and the polynomial $B(z) = A^{-1}(z) = 1 + b_1 z + b_2 z^2 + \cdots$ is bounded away from zero for $|z| \leq 1$. Therefore, the spectral density function, $f(\lambda)$, of $\{x_t\}$ satisfies $f_1 \leq f(\lambda) \leq f_2$ for some $0 < f_1 \leq f_2 < \infty$, where $-\pi < \lambda \leq \pi$. This property also ensures that $\sup_{k \geq 1} \|R(k)\| < \infty$ and $\sup_{k \geq 1} \|R^{-1}(k)\| < \infty$, where

$$R(k) = E(\mathbf{x}_n(k) \mathbf{x}'_n(k)) \tag{2.2}$$

and $\|A\|^2 = \lambda_{\max}(A'A)$ denotes the maximal eigenvalue of the matrix $A'A$. Moreover, according to Brillinger ([8], Theorem 3.8.4), (K.1)(b) and (K.1)(c) imply that $\sum_{i=1}^{\infty} |i^{1/2} b_i| < \infty$ and $\sum_{i=1}^{\infty} |i b_i| < \infty$, respectively.

The MSPE of $\hat{x}_{n+1}(k)$ can be expressed as

$$E(x_{n+1} - \hat{x}_{n+1}(k))^2 - \sigma^2 = E(\mathbf{f}(k) + \mathcal{S}(k))^2, \tag{2.3}$$

where $1 \leq k \leq K_n$,

$$\mathbf{f}(k) = \mathbf{x}'_n(k) \hat{R}_n^{-1}(k) \frac{1}{N} \sum_{j=K_n}^{n-1} \mathbf{x}_j(k) e_{j+1,k},$$

$$e_{j+1,k} = x_{j+1} + \mathbf{x}'_j(k) \mathbf{a}(k),$$

$$\mathbf{a}(k) = (a_1(k), \ldots, a_k(k))'$$

is the minimizer of $m_k(\mathbf{c}) = E(x_{k+1} + \mathbf{x}'_k(k)\mathbf{c})^2, \mathbf{c} \in R^k$, and

$$\mathcal{S}(k) = \sum_{i=1}^{\infty}(a_i - a_i(k))x_{n+1-i},$$

with $a_i(k) = 0$ for $i > k$. To simplify the notation, $\mathbf{a}(k)$ is sometimes viewed as an infinite-dimensional vector with entries $a_i(k)$, $i = 1, 2, \ldots$.

To find an asymptotic expression for $l_n(k) = E(x_{n+1} - \hat{x}_{n+1}(k))^2 - \sigma^2$, Proposition 1 below deals with the moment properties of $\hat{R}_n^{-1}(k)$, defined after (1.3). Its proof can be found in [17] [see equations (2.27) and (2.28) and Theorem 2]. For the sake of convenience, in the rest of this paper we use $C$ to denote a generic positive constant independent of the sample size $n$ and of any index with an upper (or lower) limit depending on $n$. But $C$ may depend on the distributional properties of $x_t$. It also may have different values in different places.

PROPOSITION 1. *Assume* (K.1)(a), (K.2), (K.3) *and* (K.4). *Then, for any $q > 0$,*

(2.4) $$\max_{1 \leq k \leq K_n} E\|\hat{R}_n^{-1}(k)\|^q \leq C$$

*and*

(2.5) $$\max_{1 \leq k \leq K_n} \frac{E\|\hat{R}_n^{-1}(k) - R^{-1}(k)\|^q}{(k^2/N)^{q/2}} \leq C$$

*hold for all sufficiently large $n$, where $R^{-1}(k)$ denotes the inverse of $R(k)$ [see (2.2)].*

Armed with Proposition 1, Ing and Wei ([17], Theorem 3) obtained an asymptotic expression for $l_n(k)$ which holds uniformly for all $1 \leq k \leq K_n$. This result is summarized in the following proposition.

PROPOSITION 2. *Assume* (K.1)(b), (K.2), (K.3) *and* (K.4). *Then*

(2.6) $$\lim_{n \to \infty} \max_{1 \leq k \leq K_n} \left| \frac{E(x_{n+1} - \hat{x}_{n+1}(k))^2 - \sigma^2}{L_n(k)} - 1 \right| = 0,$$

*where $L_n(k) = (k/N)\sigma^2 + \|\mathbf{a} - \mathbf{a}(k)\|_R^2$, $\mathbf{a} = (a_1, a_2, \ldots)'$, $\mathbf{a}(k)$ is now viewed as an infinite-dimensional vector, and for an infinite-dimensional vector $\mathbf{d} = (d_1, d_2, \ldots)'$,*

$$\|\mathbf{d}\|_R^2 = \sum_{i \leq i,j \leq \infty} d_i d_j \gamma_{i-j},$$

*with $\gamma_{i-j} = E(x_i x_j)$. We also note that $\|\mathbf{a} - \mathbf{a}(k)\|_R^2 = E(\mathcal{S}^2(k))$ decreases as $k$ increases.*

The following result provides an asymptotic expression for the MSPE of the least squares predictor, $\hat{y}_{n+1}(k)$, in independent-realization settings.

PROPOSITION 3. *Assume that the assumptions of Proposition 2 hold. Then*

$$\lim_{n \to \infty} \max_{1 \leq k \leq K_n} \left| \frac{E\{(y_{n+1} - \hat{y}_{n+1}(k))^2\} - \sigma^2}{L_n(k)} - 1 \right| = 0. \tag{2.7}$$

A proof of Proposition 3 can be found in Theorem 4 of [17]. Viewing (2.6) and (2.7), both types of second-order MSPEs can be uniformly approximated by the same function, $L_n(k)$, and, hence, they are asymptotically equivalent. However, this equivalence should not be taken for granted. To see this, Ing [15] recently showed that, if the underlying process is a random walk model and the assumed model is correctly specified, then

$$\lim_{n \to \infty} \frac{E(x_{n+1} - \hat{x}_{n+1}(1))^2 - \sigma^2}{E(y_{n+1} - \hat{y}_{n+1}(1))^2 - \sigma^2} \doteq \frac{2}{13.2859}.$$

Therefore, the equivalence mentioned above does not hold in this example. For stationary AR processes, Kunitomo and Yamamoto ([20], pages 946–947) also considered a comparison between same- and independent-realization MSPEs. They showed that the difference between the terms of order $1/n$ of the two types of MSPEs can be substantial when a fixed-order and underspecified AR model is used. [Note that their conclusion does not contradict that obtained from (2.6) and (2.7), because the second-order MSPE is of order $O(1)$ in the underspecified and fixed-order case.] These comparisons show that the difference between the MSPEs in two types of forecasting settings should be carefully examined in each different situation. It can be erroneous to directly assume that the results for same-realization predictions will be the same as those for the corresponding independent case without theoretical justification.

Due to more complicated probabilistic structures, this analysis becomes more difficult when the order of the predictor is selected by a data-driven method. This situation is considered in the following section.

**3. Asymptotic expressions for the MSPEs of AIC and its variants.** Values of $S_n$, AIC, FPE, $S_p$ and $C_p$ for an AR($k$) model are defined by

$$S_n(k) = (N + 2k)\hat{\sigma}_k^2,$$

$$\text{AIC}(k) = \log \hat{\sigma}_k^2 + \frac{2k}{n},$$

$$\text{FPE}(k) = \left(\frac{n+k}{n-k}\right)\hat{\sigma}_k^2,$$

$$\text{S}_p(k) = \left(1 + \frac{k}{N-k-1}\right)\hat{\hat{\sigma}}_k^2$$

and
$$C_p(k) = N\hat{\sigma}_k^2 - (N - 2k)\hat{\hat{\sigma}}_{K_n}^2,$$
respectively, where
$$\hat{\sigma}_k^2 = \frac{1}{N} \sum_{t=K_n}^{n-1} (x_{t+1} + \hat{a}_{1,n}(k)x_t + \cdots + \hat{a}_{k,n}(k)x_{t+1-k})^2$$
and
$$\hat{\hat{\sigma}}_k^2 = \left(\frac{N}{N-k}\right)\hat{\sigma}_k^2.$$
Also define
$$\hat{k}_n^S = \operatorname*{arg\,min}_{1 \le k \le K_n} S_n(k),$$
$$\hat{k}_n^A = \operatorname*{arg\,min}_{1 \le k \le K_n} \operatorname{AIC}(k),$$
$$\hat{k}_n^F = \operatorname*{arg\,min}_{1 \le k \le K_n} \operatorname{FPE}(k),$$
$$\hat{k}_n^{S_p} = \operatorname*{arg\,min}_{1 \le k \le K_n} S_p(k)$$
and
$$\hat{k}_n^C = \operatorname*{arg\,min}_{1 \le k \le K_n} C_p(k).$$

For independent-realization predictions, Shibata ([27], Section 4), assuming (K.1)(a), (K.5), $K_n = o(n^{1/2})$ and Gaussian noise, showed that

$$(3.1) \qquad \left| \frac{E\{(y_{n+1} - \hat{y}_{n+1}(\hat{k}_n))^2 | x_1, \ldots, x_n\} - \sigma^2}{L_n(k_n^*)} - 1 \right| = o_p(1),$$

where $1 \le \hat{k}_n \le K_n$ equals $\hat{k}_n^A$, $\hat{k}_n^F$ or $\hat{k}_n^S$, and $k_n^*$ is defined implicitly through $L_n(k_n^*) = \min_{1 \le k \le K_n} L_n(k)$. (Note that $k_n^* \to \infty$, provided $K_n \to \infty$ and (K.5) holds; see [27], page 154.) However, since, as mentioned in the first section, this paper focuses on the unconditional MSPE, an unconditional version of (3.1) is given in the following theorem.

THEOREM 1. *Let the assumptions of Proposition 2 and* (K.5) *hold. Then*
$$\lim_{n \to \infty} \frac{E(y_{n+1} - \hat{y}_{n+1}(\hat{k}_n))^2 - \sigma^2}{L_n(k_n^*)} = 1,$$
*where* $\hat{k}_n = \hat{k}_n^A, \hat{k}_n^F, \hat{k}_n^C, \hat{k}_n^{S_p}$ *or* $\hat{k}_n^S$.

REMARK 2. It is unclear from Shibata's [27] paper whether (3.1) holds with $\hat{k}_n = \hat{k}_n^C$ or $\hat{k}_n^{S_p}$. In fact, both $C_p$'s and $S_p$'s predictive abilities in AR($\infty$) models have seldom been discussed in the literature. On the other hand, $C_p$ and AIC have been proven to be asymptotically equivalent in the regression model with infinitely many parameters; see, for example, [28] and [25]. Under a similar situation, Breiman and Freedman [7] also established $S_p$'s asymptotic optimality for prediction (see Section 6 for more details). These previous results motivated us to include $C_p$ and $S_p$ in the analysis.

Theorem 1 shows that, for independent-realization settings, the second-order (unconditional) MSPE of the least squares predictor with the order selected by AIC, $S_n(k)$, FPE, $S_p$ or $C_p$ can ultimately achieve the best compromise between model complexity, $(k/N)\sigma^2$, and goodness of fit, $\|\mathbf{a} - \mathbf{a}(k)\|_R^2$. This result led us to ask whether AIC still possesses a similar property for same-realization predictions. Since the model selection criteria, estimated parameters and future observations are all stochastically dependent in same-realization settings, we impose the following assumption on $L_n(k)$ in order to simplify the dependent structures among these components.

(K.6) For any $\xi > 0$, there is an exponent $\theta = \theta(\xi)$ with $0 \leq \theta < 1$ such that, for all large $n$ and all $k \in A_{n,\theta} = \{k : 1 \leq k \leq K_n, |k - k_n^*| \geq k_n^{*\theta}\}$,

$$(3.2) \qquad k_n^{*\xi} \frac{N(L_n(k) - L_n(k_n^*))}{|k - k_n^*|} \geq \bar{C} > 0,$$

where $K_n$ satisfies (K.4) and $\bar{C}$ is some positive constant independent of $n$.

Note that if $\{x_t\}$ is a stationary AR model of finite order, then (3.2) holds automatically. When $\{x_t\}$ is truly a stationary AR($\infty$) model, the following two examples also show that (3.2) is flexible enough to accommodate a variety of applications.

EXAMPLE 1 (Exponential-decay case). Assume that, for all $k = 0, 1, \ldots$, the AR coefficients satisfy

$$(3.3) \qquad C_1 k^{-\theta_1} e^{-\beta k} \leq \sum_{i \geq k} a_i^2 \leq C_2 k^{\theta_1} e^{-\beta k},$$

where $\theta_1$ is some nonnegative number, and $\beta, C_2$ and $C_1$ are some positive numbers with $C_2 \geq C_1$. Note that (3.3) is satisfied by any causal and invertible ARMA($p$, $q$) model with $q > 0$. It is shown in the Appendix that (3.3) and (K.4) yield

$$k_n^* = \frac{1}{\beta} \log N + O(\log_2 N),$$

where log denotes the natural logarithm and $\log_2 N = \log(\log N)$. This result and (3.3) further ensure that, for any $0 < \eta < 1$ and all $|k - k_n^*| > k_n^{*\eta}$,

$$(3.4) \qquad \frac{N(L_n(k) - L_n(k_n^*))}{|k - k_n^*|} \geq C_3$$

holds for all sufficiently large $n$ and some positive number $C_3$ independent of $n$. Hence, for any $\xi > 0$, (3.2) is satisfied with any $0 < \theta < 1$. For a proof of (3.4), see the Appendix.

EXAMPLE 2 (Algebraic-decay case). Assume that, for all $k = 0, 1, \ldots$,

$$(3.5) \qquad (C_4 - M_1 k^{-\xi_1}) k^{-\beta} \leq \|\mathbf{a} - \mathbf{a}(k)\|_R^2 \leq (C_4 + M_1 k^{-\xi_1}) k^{-\beta},$$

where $C_4, M_1, \xi_1$ and $\beta$ are some positive numbers. Note that for independent-realization predictions, Shibata ([27], page 162) gave a similar condition,

$$(3.6) \qquad \|\mathbf{a} - \mathbf{a}(k)\|_R^2 = C_5 k^{-\beta},$$

to illustrate that AIC is strictly better than the other criteria that have different weights for penalizing the number of regressors in the model. However, since (3.6) imposes a rather restrictive limitation on $\|\mathbf{a} - \mathbf{a}(k)\|_R^2$, we use (3.5) to replace it. Under (K.4) and (3.5) with $\xi_1 \geq 2$ and $\beta > 1 + \delta_1$ [note that $\delta_1$, defined in (K.4), can be an arbitrarily small positive number], we show in the Appendix that

$$k_n^* = \left(\frac{\sigma^2}{NC_4\beta}\right)^{-1/(\beta+1)} + O(1),$$

and for some positive number $C_6$ and all $|k - k_n^*| > C_6$,

$$(3.7) \qquad \frac{N(L_n(k) - L_n(k_n^*))}{|k - k_n^*|} \geq C_7\left(\left|\frac{k - k_n^*}{k_n^*}\right| \wedge 1\right)$$

holds for all sufficiently large $n$ and some positive number $C_7$ independent of $n$, where, for real numbers $a$ and $b$, $a \wedge b = a$ if $a \leq b$ and $a \wedge b = b$ if $a > b$. Therefore, for any $\xi > 0$, (3.2) is satisfied with any $1 - \min\{\xi, 1\} < \theta < 1$.

The above discussion shows that assumption (K.6) is quite natural from both practical and theoretical points of view, since it includes the ARMA models (which are the most used short-memory time series models in practice) and the AR models with algebraic-decay coefficients (which are of much theoretical interest in the context of model selection) as special cases. Technically speaking, (3.2) gives $L_n(k)$ a basin-like shape such that, for $k$ distant from the bottom (falling into $A_{n,\theta}$), the probability of $\{\hat{k}_n = k\}$, with $\hat{k}_n = \hat{k}_n^A, \hat{k}_n^F, \hat{k}_n^C, \hat{k}_n^{S_p}$ or $\hat{k}_n^S$, is "sufficiently" small (see the proof of Theorem 2 for more details). Now the main result of this section is stated as follows.

THEOREM 2. *Let the assumptions of Theorem* 1 *and* (K.6) *hold. Then for* same-realization predictions, *one has*

$$\lim_{n \to \infty} \frac{E(x_{n+1} - \hat{x}_{n+1}(\hat{k}_n))^2 - \sigma^2}{L_n(k_n^*)} = 1, \tag{3.8}$$

*where* $\hat{k}_n = \hat{k}_n^A, \hat{k}_n^F, \hat{k}_n^C, \hat{k}_n^{S_p}$ *or* $\hat{k}_n^S$.

In the literature the penalty for the number of regressors in the model sometimes has a weight different from that used in AIC. Following [3], we now consider $\text{AIC}_\alpha(k)$, where $\alpha > 1$,

$$\text{AIC}_\alpha(k) = \log \hat{\sigma}_k^2 + \frac{\alpha k}{n},$$

$$L_n^{(\alpha)}(k) = \frac{(\alpha - 1)k\sigma^2}{N} + \|\mathbf{a} - \mathbf{a}(k)\|_R^2,$$

$$\hat{k}_n^{A_\alpha} = \underset{1 \le k \le K_n}{\arg\min} \, \text{AIC}_\alpha(k)$$

and

$$k_n^{*(\alpha)} = \underset{1 \le k \le K_n}{\arg\min} \, L_n^{(\alpha)}(k).$$

To investigate the performances of $\text{AIC}_\alpha(k)$, $\alpha > 1$, for same-realization predictions, we need the following analogy of (K.6).

(K.6′) For any $\xi > 0$, there exists an exponent $\theta = \theta(\xi)$ with $0 \le \theta < 1$ such that, for all large $n$ and all $k \in A_{n,\theta}^{(\alpha)} = \{k : 1 \le k \le K_n, |k - k_n^{*(\alpha)}| \ge (k_n^{*(\alpha)})^\theta\}$,

$$(k_n^{*(\alpha)})^\xi \frac{N(L_n^{(\alpha)}(k) - L_n^{(\alpha)}(k_n^{*(\alpha)}))}{|k - k_n^{*(\alpha)}|} \ge \bar{C} > 0, \tag{3.9}$$

where $\alpha > 1$, $K_n$ satisfies (K.4) and $\bar{C}$ is some positive number independent of $n$.

For any $\alpha > 1$, it is easy to see that (3.9) is fulfilled by finite-order stationary AR models. By arguments similar to those given in the Appendix, we can also show that, for any $\alpha > 1$, (3.9) is satisfied by stationary AR($\infty$) models with coefficients which obey (3.3) or (3.5) (with $\beta > 1 + \delta_1$ and $\xi_1 \ge 2$). Therefore, assumption (K.6′), like assumption (K.6), is reasonable for a wide range of applications. In Corollary 1 we obtain an asymptotic expression for the MSPE of $\hat{x}_{n+1}(\hat{k}_n^{A_\alpha})$, $\alpha > 1$, in same-realization settings.

COROLLARY 1. *Let the assumptions of Theorem 2 hold with* (K.6) *replaced by* (K.6′). *Then*

$$\lim_{n\to\infty} \frac{E(x_{n+1} - \hat{x}_{n+1}(\hat{k}_n^{A_\alpha}))^2 - \sigma^2}{L_n(k_n^{*(\alpha)})} = 1. \tag{3.10}$$

REMARK 3. By arguments like those used to prove Corollary 1, (3.10) still holds with $\hat{k}_n^{A_\alpha}$ replaced by $\hat{k}_n^{F_\alpha}$ or $\hat{k}_n^{(\alpha)}$, where $\alpha > 1$,

$$\hat{k}_n^{F_\alpha} = \underset{1 \leq k \leq K_n}{\arg\min}\, \mathrm{FPE}_\alpha(k) = \underset{1 \leq k \leq K_n}{\arg\min} \left(1 + \frac{\alpha k}{n}\right)\hat{\sigma}_k^2$$

and

$$\hat{k}_n^{(\alpha)} = \underset{1 \leq k \leq K_n}{\arg\min}\, S_n^{(\alpha)}(k) = \underset{1 \leq k \leq K_n}{\arg\min} (N + \alpha k)\hat{\sigma}_k^2.$$

Therefore, $\mathrm{AIC}_\alpha(k)$, $\mathrm{FPE}_\alpha(k)$ and $S_n^{(\alpha)}(k)$ are (asymptotically) equally efficient for the same choice of $\alpha$. Note that $\mathrm{FPE}_\alpha(k)$ and $S_n^{(\alpha)}(k)$ were first introduced by Bhansali and Downham [6] and Shibata [27], respectively.

To illustrate Corollary 1, we first consider a special case of (3.3),

$$C_1 e^{-\beta k} \leq \sum_{i \geq k} a_i^2 \leq C_2 e^{-\beta k}, \tag{3.11}$$

where $\beta, C_1$ and $C_2$ are some positive numbers with $C_2 \geq C_1$. Condition (3.11) is satisfied by any causal and invertible ARMA$(p, q)$ process with $q > 0$. Under (3.11), it can be shown that, for any $\alpha > 1$,

$$\frac{L_n(k_n^{*(\alpha)})}{L_n(k_n^*)} = 1.$$

This fact and (3.10) yield that, for any two positive numbers $\alpha_1$ and $\alpha_2$ larger than 1, $\mathrm{AIC}_{\alpha_1}(k)$ and $\mathrm{AIC}_{\alpha_2}(k)$ are asymptotically equivalent, namely,

$$\lim_{n\to\infty} \frac{E(x_{n+1} - \hat{x}_{n+1}(\hat{k}_n^{A_{\alpha_2}}))^2 - \sigma^2}{E(x_{n+1} - \hat{x}_{n+1}(\hat{k}_n^{A_{\alpha_1}}))^2 - \sigma^2} = 1. \tag{3.12}$$

Next consider the algebraic-decay case (3.5) with $\xi_1 \geq 2$ and $\beta > 1 + \delta_1$. By arguments similar to those used for obtaining (3.7) and Case I of [27], page 162, one has, for $1 < \alpha_2 < \alpha_1 \leq 2$ or $2 \leq \alpha_1 < \alpha_2 < \infty$,

$$\liminf_{n\to\infty} \frac{L_n(k_n^{*(\alpha_2)})}{L_n(k_n^{*(\alpha_1)})} > 1$$

and, hence, for $1 < \alpha_2 < \alpha_1 \leq 2$ or $2 \leq \alpha_1 < \alpha_2 < \infty$,

$$(3.13) \qquad \liminf_{n \to \infty} \frac{E(x_{n+1} - \hat{x}_{n+1}(\hat{k}_n^{A_{\alpha_2}}))^2 - \sigma^2}{E(x_{n+1} - \hat{x}_{n+1}(\hat{k}_n^{A_{\alpha_1}}))^2 - \sigma^2} > 1.$$

Inequality (3.13) and Corollary 1 together imply that AIC *asymptotically dominates* $\text{AIC}_\alpha$, $\alpha \neq 2$, in the sense that

$$(3.14) \qquad \liminf_{n \to \infty} \frac{E(x_{n+1} - \hat{x}_{n+1}(\hat{k}_n^{A_\alpha}))^2 - \sigma^2}{E(x_{n+1} - \hat{x}_{n+1}(\hat{k}_n^{A}))^2 - \sigma^2} \geq 1,$$

with strict inequality holding for at least the algebraic-decay case (3.5).

Before leaving this section, we note that (3.12)–(3.14) seem to be the first results that can evaluate (compare) the performances of $\text{AIC}_\alpha(k)$, $\alpha > 1$, from the same-realization prediction point of view. For more applications of these results, see Section 4. In addition to (3.5), we have also found a similar case, $\|\mathbf{a} - \mathbf{a}(k)\|_R^2 = C_5(\log k)^{\theta_3} k^{-\beta}$, with $C_5 > 0$, $-\infty < \theta_3 < \infty$, and $\beta > 1 + \delta_1$, where (3.14) holds with strict inequality only. However, to gain a deeper understanding of AIC it would be interesting to identify more $\text{AR}(\infty)$ models which can lead to the same property.

**4. Performances of AIC and its variants for independent- and same-realization predictions.** Based on the results obtained in Section 3, this section aims to investigate how well AIC (or its variants) works for independent- and same-realization predictions. Let $\hat{k}_n \in \{1, 2, \ldots, K_n\}$ be determined by a certain order selection criterion with $K_n$ satisfying assumption (K.4). Define

$$PE(\hat{x}_{n+1}(\hat{k}_n)) = \frac{E(x_{n+1} - \hat{x}_{n+1}(\hat{k}_n))^2 - \sigma^2}{\min_{1 \leq k \leq K_n} E(x_{n+1} - \hat{x}_{n+1}(k))^2 - \sigma^2}.$$

We say that $\hat{k}_n$ is asymptotically efficient for same-realization predictions if

$$(4.1) \qquad \limsup_{n \to \infty} PE(\hat{x}_{n+1}(\hat{k}_n)) \leq 1.$$

Similarly, for independent-realization predictions, define $PEI(\hat{y}_{n+1}(\hat{k}_n))$ as

$$PEI(\hat{y}_{n+1}(\hat{k}_n)) = \frac{E(y_{n+1} - \hat{y}_{n+1}(\hat{k}_n))^2 - \sigma^2}{\min_{1 \leq k \leq K_n} E(y_{n+1} - \hat{y}_{n+1}(k))^2 - \sigma^2}.$$

We say that $\hat{k}_n$ is asymptotically efficient for independent-realization predictions if

$$(4.2) \qquad \limsup_{n \to \infty} PEI(\hat{y}_{n+1}(\hat{k}_n)) \leq 1.$$

Inequality (4.1) [(4.2)] says that, if $\hat{k}_n$ is determined by an asymptotically efficient criterion, then the relative prediction efficiency of the best predictor (from the MSPE point of view) among $\{\hat{x}_{n+1}(1),\ldots,\hat{x}_{n+1}(K_n)\}$ $[\{\hat{y}_{n+1}(1),\ldots,\hat{y}_{n+1}(K_n)\}]$ of $\hat{x}_{n+1}(\hat{k}_n)$ $[\hat{y}_{n+1}(\hat{k}_n)]$, that is, $PE(\hat{x}_{n+1}(\hat{k}_n))$ $(PEI(\hat{y}_{n+1}(\hat{k}_n)))$, will ultimately not exceed 1.

Proposition 3 and Theorem 1 yield that

$$\lim_{n\to\infty} PEI(\hat{y}_{n+1}(\hat{k}_n)) = 1, \tag{4.3}$$

where $\hat{k}_n = \hat{k}_n^A, \hat{k}_n^F, \hat{k}_n^S, \hat{k}_n^{S_p}$ or $\hat{k}_n^C$. Therefore, AIC, FPE, $S_n(k)$, $S_p$ and $C_p$ are all asymptotically efficient for independent-realization predictions. According to Proposition 2 and Theorem 2, we have

$$\lim_{n\to\infty} PE(\hat{x}_{n+1}(\hat{k}_n)) = 1, \tag{4.4}$$

where $\hat{k}_n = \hat{k}_n^A, \hat{k}_n^F, \hat{k}_n^S, \hat{k}_n^{S_p}$ or $\hat{k}_n^C$, which shows that these criteria are also asymptotically efficient for same-realization predictions. In addition, if it is already known that the exponential-decay case (3.11) holds, then (3.12) and Remark 3 suggest that more options, for example, $\text{AIC}_\alpha$ [$\text{FPE}_\alpha$, $S_n^{(\alpha)}(k)$], with any $\alpha > 1$, are available for achieving asymptotic efficiency. However, AIC and its variants cannot be replaced by $\text{AIC}_\alpha$ [$\text{FPE}_\alpha$, $S_n^{(\alpha)}(k)$], $\alpha \neq 2$, in general, since (3.13) and (4.4) imply that the latter criterion is not asymptotically efficient in the algebraic-decay case.

To gain further insight into the practical implications of asymptotically efficient criteria for same-realization predictions, a simulation study is conducted. Let observations be generated from an ARMA(1, 1) model

$$x_{t+1} = \phi_0 x_t + \varepsilon_t - \theta_0 \varepsilon_{t-1},$$

where the $\varepsilon_t$'s are independent and identically $\mathcal{N}(0,1)$ distributed, $\phi_0 = \pm 0.9, \pm 0.7, \pm 0.5$ and $\theta_0 = \pm 0.8, \pm 0.6$. For each combination $(\phi_0, \theta_0)$, the empirical estimates of $PE(\hat{x}_{n+1}(\hat{k}_n^A))$, denoted by $\widehat{PE}(\hat{x}_{n+1}(\hat{k}_n^A))$, are obtained based on 20,000 replications for $(n, K_n) = (60, 7), (120, 10), (200, 14), (500, 22)$ and $(1000, 31)$. (Note that $K_n$ here is set to the largest integer $\leq n^{1/2}$.) In addition, the empirical estimates of

$$\gamma_{\text{opt}}(n, K_n) = \frac{\min_{1 \leq k \leq K_n} E(x_{n+1} - \hat{x}_{n+1}(k))^2 - \sigma^2}{\min_{1 \leq k \leq 6} E(x_{61} - \hat{x}_{61}(k))^2 - \sigma^2},$$

denoted by $\hat{\gamma}_{\text{opt}}(n, K_n)$, with $\sigma^2 = 1$ and $(n, K_n) = (60, 7), (120, 10), (200, 14), (500, 22)$ and $(1000, 31)$, are also obtained based on the same ARMA(1, 1) model and 20,000 replications. [Note that $\gamma_{\text{opt}}(n, K_n)$ is used to illustrate how fast $\min_{1 \leq k \leq K_n} E(x_{n+1} - \hat{x}_{n+1}(k))^2 - \sigma^2$ decreases as $n$ and $K_n$ simultaneously increase.] According to the rate of convergence of $\hat{\gamma}_{\text{opt}}(n, K_n)$, these empirical results (which are summarized in Table 1) can be classified into three categories. [Since $\gamma_{\text{opt}}(60, 7) = 1$, $\hat{\gamma}_{\text{opt}}(60, 7)$ is set to 1 in Table 1.]

We first observe that the fast rate of decrease of $\hat{\gamma}_{\text{opt}}(n, K_n)$ clearly occurs when $\text{sgn}(\phi_0) \neq \text{sgn}(\theta_0)$, $(\phi_0, \theta_0) = (0.9, 0.6)$, or $(\phi_0, \theta_0) = (-0.9, -0.6)$, where, for a nonzero real number $a$, $\text{sgn}(a) = 1$ if $a > 0$ and $\text{sgn}(a) = -1$ if $a < 0$. In these cases, we also observe that the rate of convergence of $\widehat{PE}(\hat{x}_{n+1}(\hat{k}_n^A))$ is very slow and the values fluctuate around a certain number which is not distant from 1. In particular, if $\theta_0 = \pm 0.8$, the values fluctuate around 1.25, while, if $\theta_0 = \pm 0.6$, the values fluctuate around 1.5 (or slightly higher). The second category contains those parameter combinations satisfying $|\phi_0 - \theta_0| = 0.1$. The decreasing rate of $\hat{\gamma}_{\text{opt}}(n, K_n)$ in this category is obviously slower than that in the first category. On the other hand, except in the cases where $(\phi_0, \theta_0) = (0.9, 0.8)$ and $(-0.9, -0.8)$, the advantage of increasing $n$ and $K_n$ in reducing the values of $\widehat{PE}(\hat{x}_{n+1}(\hat{k}_n^A))$ becomes quite significant and the ratio $\widehat{PE}(\hat{x}_{1001}(\hat{k}_{1000}^A))/\widehat{PE}(\hat{x}_{61}(\hat{k}_{60}^A))$ is smaller than $2/3$. When $(\phi_0, \theta_0) = (0.9, 0.8)$ and $(-0.9, -0.8)$, the values of $\widehat{PE}(\hat{x}_{n+1}(\hat{k}_n^A))$ are smaller than those in the other cases in this category. However, the reduction in the values of $\widehat{PE}(\hat{x}_{n+1}(\hat{k}_n^A))$ is also much smaller (only a slightly decreasing trend can be observed). Another observation regarding this category is that, as $(n, K_n)$ increases to $(1000, 31)$, $\widehat{PE}(\hat{x}_{n+1}(\hat{k}_n^A))$ decreases to a value around 1.5 if $\theta_0 = \pm 0.8$, and decreases to a value around 1.95 if $\theta_0 = \pm 0.6$. The third category contains the remaining parameter combinations, namely, $(\phi_0, \theta_0) = (0.5, 0.8)$ and $(-0.5, -0.8)$. The rate of convergence of $\hat{\gamma}_{\text{opt}}(n, K_n)$ in this category is intermediate (slower than in the first category, but faster than in the second category), while the rate for $\widehat{PE}(\hat{x}_{n+1}(\hat{k}_n^A))$ is slow and its value fluctuates around 1.5. In summary, although the convergence rate of $\widehat{PE}(\hat{x}_{n+1}(\hat{k}_n^A))$ is slow when $|\phi_0 - \theta_0| \geq 0.3$ (which includes the first and the third categories), the relatively small values of $\widehat{PE}(\hat{x}_{n+1}(\hat{k}_n^A))$, accompanied by a (very) fast rate of decrease of $\hat{\gamma}_{\text{opt}}(n, K_n)$, show that, even in finite sample situations, AIC also yields a satisfactory same-realization prediction for $AR(\infty)$ models, provided the number of candidate models is allowed to increase with the sample size. On the other hand, when $|\phi_0 - \theta_0| = 0.1$, the (prediction) efficiency of AIC is not satisfactory (except in the cases where $|\phi_0| = 0.9$) and the reduction in the values of $\hat{\gamma}_{\text{opt}}(n, K_n)$ through increasing $n$ and $K_n$ also becomes relatively insignificant. However, the efficiency of AIC can be substantially improved (except in the cases where $|\phi_0| = 0.9$) if $n$ and $K_n$ are allowed to increase simultaneously.

In the rest of this section, the question of whether AIC is asymptotically optimal among all order selection criteria from the MSPE point of view is investigated. This question led us to define a stronger version of asymptotic efficiency. An order selection criterion $\hat{k}_n^{\text{SA}}$ with $1 \leq \hat{k}_n^{\text{SA}} \leq K_n$ is said to be *strongly* asymptotically efficient for same-realization predictions if

$$(4.5) \qquad \lim_{n \to \infty} \frac{E(x_{n+1} - \hat{x}_{n+1}(\hat{k}_n^{\text{SA}}))^2 - \sigma^2}{\inf_{\hat{I}_n \in \mathcal{J}_n} E(x_{n+1} - \hat{x}_{n+1}(\hat{I}_n))^2 - \sigma^2} = 1,$$

TABLE 1
*Empirical estimates of $PE(\hat{x}_{n+1}(\hat{k}_n^A))$ and $\gamma_{\mathrm{opt}}(n, K_n)$*

| | | \multicolumn{8}{c}{$\theta_0$} |
|---|---|---|---|---|---|---|---|---|---|
| | | **0.8** | | **0.6** | | **−0.6** | | **−0.8** | |
| $\phi_0$ | $n/K_n$ | P | R | P | R | P | R | P | R |
| −0.9 | 60/7 | 1.23 | 1 | 1.46 | 1 | 1.64 | 1 | 1.68 | 1 |
| | 120/10 | 1.21 | 0.56 | 1.46 | 0.54 | 1.58 | 0.58 | 1.55 | 0.71 |
| | 200/14 | 1.29 | 0.36 | 1.56 | 0.34 | 1.67 | 0.38 | 1.52 | 0.54 |
| | 500/22 | 1.25 | 0.17 | 1.56 | 0.15 | 1.63 | 0.17 | 1.53 | 0.34 |
| | 1000/31 | 1.29 | 0.09 | 1.54 | 0.08 | 1.58 | 0.10 | 1.46 | 0.17 |
| −0.7 | 60/7 | 1.23 | 1 | 1.44 | 1 | 2.90 | 1 | 2.49 | 1 |
| | 120/10 | 1.25 | 0.57 | 1.49 | 0.54 | 2.48 | 0.68 | 2.37 | 0.70 |
| | 200/14 | 1.29 | 0.37 | 1.61 | 0.35 | 2.25 | 0.51 | 1.95 | 0.59 |
| | 500/22 | 1.31 | 0.17 | 1.53 | 0.16 | 1.97 | 0.27 | 1.65 | 0.38 |
| | 1000/31 | 1.28 | 0.10 | 1.56 | 0.09 | 1.93 | 0.17 | 1.62 | 0.23 |
| −0.5 | 60/7 | 1.23 | 1 | 1.50 | 1 | 3.48 | 1 | 1.49 | 1 |
| | 120/10 | 1.25 | 0.57 | 1.62 | 0.53 | 3.01 | 0.65 | 1.47 | 0.67 |
| | 200/14 | 1.33 | 0.36 | 1.57 | 0.35 | 2.78 | 0.47 | 1.53 | 0.46 |
| | 500/22 | 1.33 | 0.17 | 1.56 | 0.16 | 2.29 | 0.27 | 1.44 | 0.21 |
| | 1000/31 | 1.26 | 0.10 | 1.56 | 0.09 | 1.99 | 0.17 | 1.42 | 0.13 |
| 0.5 | 60/7 | 1.55 | 1 | 3.10 | 1 | 1.49 | 1 | 1.25 | 1 |
| | 120/10 | 1.51 | 0.67 | 2.98 | 0.63 | 1.59 | 0.55 | 1.23 | 0.56 |
| | 200/14 | 1.48 | 0.46 | 2.86 | 0.45 | 1.55 | 0.38 | 1.31 | 0.37 |
| | 500/22 | 1.47 | 0.22 | 2.45 | 0.26 | 1.61 | 0.16 | 1.28 | 0.17 |
| | 1000/31 | 1.41 | 0.13 | 1.99 | 0.16 | 1.57 | 0.09 | 1.32 | 0.10 |
| 0.7 | 60/7 | 2.71 | 1 | 2.97 | 1 | 1.55 | 1 | 1.25 | 1 |
| | 120/10 | 2.31 | 0.71 | 2.56 | 0.62 | 1.58 | 0.53 | 1.25 | 0.56 |
| | 200/14 | 1.92 | 0.62 | 2.31 | 0.48 | 1.61 | 0.36 | 1.29 | 0.37 |
| | 500/22 | 1.79 | 0.37 | 2.04 | 0.27 | 1.53 | 0.16 | 1.28 | 0.18 |
| | 1000/31 | 1.56 | 0.24 | 1.95 | 0.16 | 1.44 | 0.09 | 1.31 | 0.10 |
| 0.9 | 60/7 | 1.75 | 1 | 1.58 | 1 | 1.43 | 1 | 1.24 | 1 |
| | 120/10 | 1.56 | 0.70 | 1.61 | 0.57 | 1.50 | 0.53 | 1.23 | 0.57 |
| | 200/14 | 1.57 | 0.51 | 1.66 | 0.37 | 1.58 | 0.32 | 1.31 | 0.37 |
| | 500/22 | 1.49 | 0.29 | 1.68 | 0.17 | 1.54 | 0.15 | 1.31 | 0.17 |
| | 1000/31 | 1.48 | 0.17 | 1.57 | 0.10 | 1.47 | 0.08 | 1.29 | 0.10 |

NOTE. Column P denotes the empirical estimates of $PE(\hat{x}_{n+1}(\hat{k}_n^A))$ and column R denotes the empirical estimates of $\gamma_{\mathrm{opt}}(n, K_n)$.

and is said to be *strongly* asymptotically efficient for independent-realization predictions if

$$\text{(4.6)} \qquad \lim_{n \to \infty} \frac{E(y_{n+1} - \hat{y}_{n+1}(\hat{k}_n^{\mathrm{SA}}))^2 - \sigma^2}{\inf_{\hat{I}_n \in \mathcal{J}_n} E(y_{n+1} - \hat{y}_{n+1}(\hat{I}_n))^2 - \sigma^2} = 1.$$

Here $\mathcal{J}_n$ in (4.5) and (4.6) is the family of all $\mathcal{G}_n$-measurable random variables taking values on $\{1, 2, \ldots, K_n\}$, and $\mathcal{G}_n$ is the $\sigma$-algebra generated by $\{x_1, \ldots, x_n\}$. Observe that, for any random variable $\hat{I}_n \in \mathcal{J}_n$,

$$E(x_{n+1} - \hat{x}_{n+1}(\hat{I}_n))^2 = E\left\{\sum_{k=1}^{K_n} E[(x_{n+1} - \hat{x}_{n+1}(k))^2 | x_1, \ldots, x_n] I_{\{\hat{I}_n = k\}}\right\}$$

$$\geq E\left\{\min_{1 \leq k \leq K_n} E[(x_{n+1} - \hat{x}_{n+1}(k))^2 | x_1, \ldots, x_n]\right\},$$

where $I_{\{\hat{I}_n = k\}} = 1$ if $\hat{I}_n = k$ and $I_{\{\hat{I}_n = k\}} = 0$ if $\hat{I}_n \neq k$. Also notice that the minimizer of

$$E[(x_{n+1} - \hat{x}_{n+1}(k))^2 | x_1, \ldots, x_n], \qquad k = 1, 2, \ldots, K_n,$$

is a member of $\mathcal{J}_n$. Therefore

$$\inf_{\hat{I}_n \in \mathcal{J}_n} E(x_{n+1} - \hat{x}_{n+1}(\hat{I}_n))^2 = E\left\{\min_{1 \leq k \leq K_n} E[(x_{n+1} - \hat{x}_{n+1}(k))^2 | x_1, \ldots, x_n]\right\}$$

and, hence, (4.5) can be rewritten as

$$(4.7) \quad \lim_{n \to \infty} \frac{E(x_{n+1} - \hat{x}_{n+1}(\hat{k}_n^{\text{SA}}))^2 - \sigma^2}{E\{\min_{1 \leq k \leq K_n} E[(x_{n+1} - \hat{x}_{n+1}(k))^2 | x_1, \ldots, x_n] - \sigma^2\}} = 1.$$

Similarly, (4.6) can be rewritten as

$$(4.8) \quad \lim_{n \to \infty} \frac{E(y_{n+1} - \hat{y}_{n+1}(\hat{k}_n^{\text{SA}}))^2 - \sigma^2}{E\{\min_{1 \leq k \leq K_n} E[(y_{n+1} - \hat{y}_{n+1}(k))^2 | x_1, \ldots, x_n] - \sigma^2\}} = 1.$$

To examine whether AIC satisfies (4.5) [or (4.7)], it suffices to check whether

$$r_n^* = \frac{E(x_{n+1} - \hat{x}_{n+1}(\hat{k}_n^A))^2 - \sigma^2}{E\{\min_{1 \leq k \leq K_n} E[(x_{n+1} - \hat{x}_{n+1}(k))^2 | x_1, \ldots, x_n] - \sigma^2\}}$$

TABLE 2
*Simulation results for $r_n^*$*

| | | $\theta_0$ | | | |
|---|---|---|---|---|---|
| $n$ | $K_n$ | 0.8 | 0.6 | $-0.6$ | $-0.8$ |
| 60 | 7 | 5.61 | 6.04 | 6.48 | 5.60 |
| 120 | 10 | 7.44 | 7.82 | 7.97 | 8.16 |
| 200 | 14 | 9.68 | 9.17 | 9.66 | 9.17 |
| 500 | 22 | 13.10 | 12.30 | 12.48 | 13.54 |
| 1000 | 31 | 16.18 | 13.56 | 13.84 | 15.72 |

converges to 1. Instead of evaluating $r_n^*$ theoretically, empirical estimates of $r_n^*$, denoted by $\hat{r}_n^*$, with $(n, K_n) = (60, 7), (120, 10), (200, 14), (500, 22), (1000, 31)$, are obtained based on the MA(1) model

$$x_t = \epsilon_t - \theta_0 \epsilon_{t-1} \tag{4.9}$$

and 20,000 replications. Here we take the noise $\{\epsilon_t\}$ to be i.i.d. $\mathcal{N}(0,1)$ and use the parameter values $\theta_0 = 0.8, 0.6, -0.6, -0.8$. These estimates are summarized in Table 2. It can be seen from this table that when $(n, K_n) = (60, 7)$, values of $\hat{r}_n^*$ are larger than 5.5 for all $\theta_0$. Moreover, $\hat{r}_n^*$ increases considerably as $n$ and $K_n$ grow. In particular, when $n$ and $K_n$ increase to 1000 and 31, respectively, all values of $\hat{r}_n^*$ are larger than 13.5. Viewing the relatively moderate values given in Table 1, Table 2 suggests that it seems very difficult for AIC to achieve (4.7). This is a somewhat *different* situation from that encountered in independent-realization settings. To see this, let us restrict ourselves to model (4.9) again. Motivated by (4.8), the empirical estimates of

$$\begin{aligned} r_{I,n}^* &= \frac{E(y_{n+1} - \hat{y}_{n+1}(\hat{k}_n^A))^2 - \sigma^2}{\inf_{\hat{I}_n \in \mathcal{J}_n} E(y_{n+1} - \hat{y}_{n+1}(\hat{I}_n))^2 - \sigma^2} \\ &= \frac{E(y_{n+1} - \hat{y}_{n+1}(\hat{k}_n^A))^2 - \sigma^2}{E\{\min_{1 \leq k \leq K_n} E[(y_{n+1} - \hat{y}_{n+1}(k))^2 | x_1, \ldots, x_n] - \sigma^2\}}, \end{aligned}$$

denoted by $\hat{r}_{I,n}^*$, with $(n, K_n) = (60, 7), (120, 10), (200, 14), (500, 22), (1000, 31)$ and $\theta_0 = 0.8, 0.6, -0.6, -0.8$, are obtained based on 20,000 replications (see Table 3). Table 3 shows that values of $\hat{r}_{I,n}^*$, like values of $\widehat{PE}(\hat{x}_{n+1}(\hat{k}_n^A))$ in Table 1, are not distant from 1, particularly for large $n$, $K_n$ and $|\theta_0|$. In fact, by (5.36) (see Section 5) and Theorem 1, it can be shown that

$$\lim_{n \to \infty} r_{I,n}^* = 1, \tag{4.10}$$

provided the assumptions of Theorem 1 hold. Therefore, for independent-realization predictions, AIC is asymptotically efficient, as well as strongly asymptotically efficient. [Note that (4.10) also holds with $\hat{k}_n^A$ replaced by $\hat{k}_n^S$, $\hat{k}_n^F$, $\hat{k}_n^{S_p}$ or $\hat{k}_n^C$.] Related to (4.10), but from a conditional MSPE point of view, Shibata [27] showed that, for a Gaussian AR($\infty$) process,

$$\frac{E\{(y_{n+1} - \hat{y}_{n+1}(\hat{k}_n))^2 | x_1, \ldots, x_n\} - \sigma^2}{\inf_{\hat{I}_n \in \mathcal{J}_n} E[(y_{n+1} - \hat{y}_{n+1}(\hat{I}_n))^2 | x_1, \ldots, x_n] - \sigma^2} - 1 = o_p(1) \tag{4.11}$$

holds with $\hat{\hat{k}}_n = \hat{k}_n^A, \hat{k}_n^S$ or $\hat{k}_n^F$. An order selection criterion $\hat{\hat{k}}_n$ is said to be asymptotically efficient in Shibata's paper if it satisfies (4.11). For an equivalent definition of (4.11), see (5.37).

**5. Proofs.** We first introduce two frequently used results, Lemmas 1 and 2. The proofs of Lemmas 1 and 2 are similar to those of Lemmas 3 and 4 of [17], respectively. To save space, we skip the details.

LEMMA 1. *Assume that* (K.1)(a) *holds and* $\sup_{-\infty<t<\infty} E(|e_t|^{2q}) < \infty$ *for some* $q \geq 2$. *Let* $\{m_{i,n}\}$, $i = 0, 1, 2$, *be sequences of positive integers satisfying* $m_{2,n} \geq m_{1,n} \geq m_{0,n}$ *for all* $n \geq 1$. *Then, for all* $1 \leq k \leq m_{0,n}$,

$$(5.1) \quad E\left\|\frac{1}{\sqrt{m_n}}\sum_{j=m_{1,n}}^{m_{2,n}} \mathbf{x}_j(k)(e_{j+1,k} - e_{j+1})\right\|^q \leq Ck^{q/2}\|\mathbf{a} - \mathbf{a}(k)\|_R^q,$$

*where* $m_n = m_{2,n} - m_{1,n} + 1$, $e_{j+1,k}$ *is defined after* (2.3), $\|\mathbf{a} - \mathbf{a}(k)\|_R^2$ *is defined in Proposition* 2, *and for a* $k$-*dimensional vector* $\mathbf{v} = (v_1, \ldots, v_k)'$, $\|\mathbf{v}\|^2 = \sum_{i=1}^k v_i^2$. [*Note that if* (K.5) *holds, then* $\|\mathbf{a} - \mathbf{a}(k)\|_R^2 > 0$ *for all* $k = 1, 2, \ldots$. *In this case* (5.1) *can be expressed as*

$$\max_{1 \leq k \leq m_{0,n}} \frac{E\|(1/\sqrt{m_n})\sum_{j=m_{1,n}}^{m_{2,n}} \mathbf{x}_j(k)(e_{j+1,k} - e_{j+1})\|^q}{k^{q/2}\|\mathbf{a} - \mathbf{a}(k)\|_R^q} \leq C.]$$

LEMMA 2. *Assume that* (K.1)(a) *holds and* $\sup_{-\infty<t<\infty} E\{|e_t|^q\} < \infty$ *for* $q \geq 2$. *Let* $\{m_{i,n}\}$, $i = 0, 1, 2$, *and* $\{m_n\}$ *be defined as in Lemma* 1. *Then*

$$(5.2) \quad \max_{1 \leq k \leq m_{0,n}} (k^{-q/2}) E\left\|\frac{1}{\sqrt{m_n}}\sum_{j=m_{1,n}}^{m_{2,n}} \mathbf{x}_j(k)e_{j+1}\right\|^q \leq C$$

*and*

$$(5.3) \quad \max_{1 \leq k_1 < k_2 \leq m_{0,n}} (k_2 - k_1)^{-q/2} E\left\|\frac{1}{\sqrt{m_n}}\sum_{j=m_{1,n}}^{m_{2,n}} (\mathbf{x}_j(k_2) - \mathbf{x}_j(k_1))e_{j+1}\right\|^q \leq C,$$

*where* $\mathbf{x}_j(k_1)$ *in* (5.3) *is regarded as a* $k_2$-*dimensional vector with undefined entries set to* 0.

TABLE 3
*Simulation results for* $r_{I,n}^*$

|      |       | $\theta_0$ |      |      |      |
|------|-------|------|------|------|------|
| $n$  | $K_n$ | 0.8  | 0.6  | $-0.6$ | $-0.8$ |
| 60   | 7     | 1.54 | 2.10 | 2.12 | 1.56 |
| 120  | 10    | 1.50 | 2.05 | 2.12 | 1.49 |
| 200  | 14    | 1.57 | 1.93 | 1.90 | 1.59 |
| 500  | 22    | 1.47 | 1.79 | 1.91 | 1.44 |
| 1000 | 31    | 1.42 | 1.81 | 1.75 | 1.41 |

Lemmas 3–6, listed below, are essential tools for verifying Theorems 1 and 2. To provide motivation for Lemma 3, we note that, under (K.1)(a) and the Gaussian assumption on $\{e_t\}$, Shibata ([27], Lemmas 3.2 and 3.4) showed that, for $q = 2, 4$ and $1 \leq k \leq K_n$,

$$(5.4) \quad E\left\|\left\|\frac{1}{\sqrt{N}}\sum_{j=K_n}^{n-1}\mathbf{x}_j(k)e_{j+1}\right\|_{R^{-1}(k)}^2 - k\sigma^2\right|^q = \sum_{i=1}^{q/2} C_{i,q}k^i + O(1/N)k^q,$$

where, for a $k \times k$ symmetric matrix $A$ and a $k$-dimensional vector $\mathbf{y}$, $\|\mathbf{y}\|_A^2 = \mathbf{y}'A\mathbf{y}$, and the $C_{i,q}$'s, independent of $n$ and $k$, are some positive numbers. Under (K.1)(c) with the i.i.d. assumption on $\{e_t\}$ and the assumption that $E|e_1|^{16} < \infty$, Karagrigoriou ([19], Lemma 3.1) also gave a result similar to Shibata's. However, since they needed to calculate $e_t$'s (or $x_t$'s) $4q$th cross moments, extensions of their approaches to large $q$ cases are not easy due to heavy computational burdens. For example, in order to verify (5.4) with $q = 6$ through their approaches, one must deal with $e_t$'s (or $x_t$'s) 24th cross moments. In addition, even if (5.4) holds for large $q$'s, the remainder term, $O(1/N)k^q$, may dominate the main term, $\sum_{i=1}^{q/2} C_{i,q}k^i$ (which is of order $k^{q/2}$), in situations where $k^{q/2}/N$ is large as well. This causes another difficulty since bounding the left-hand side of (5.4) by $Ck^{q/2}$ for sufficiently large $q$ and all $1 \leq k \leq K_n$ seems indispensable for our analysis, especially for proving Theorem 2. Under a slightly stronger assumption (than Shibata's) on AR coefficients, the difficulties mentioned above are resolved by Lemma 3.

LEMMA 3. *Assume* (K.1)(b), $\sup_{-\infty < t < \infty} E|e_t|^{2q} < \infty$ *for some* $q \geq 2$, *and* $K_n = O(n^{1/2})$. *Then*

$$(5.5) \quad \max_{1 \leq k \leq K_n} k^{-q/2} E\left\|\left\|\frac{1}{\sqrt{N}}\sum_{j=K_n}^{n-1}\mathbf{x}_j(k)e_{j+1}\right\|_{R^{-1}(k)}^2 - k\sigma^2\right|^q \leq C.$$

PROOF. First observe that

$$E\left\{\left|k^{-1/2}\left(\left\|\frac{1}{\sqrt{N}}\sum_{j=K_n}^{n-1}\mathbf{x}_j(k)e_{j+1}\right\|_{R^{-1}(k)}^2 - k\sigma^2\right)\right|^q\right\}$$

$$= E\left\{\left|\frac{1}{Nk^{1/2}}\sum_{j=K_n}^{n-1}(\mathbf{x}'_j(k)R^{-1}(k)\mathbf{x}_j(k)e_{j+1}^2 - k\sigma^2)\right.\right.$$

$$\left.\left. + \frac{2}{Nk^{1/2}}\sum_{l=K_n+1}^{n-1}\sum_{j=K_n}^{l-1}\mathbf{x}'_j(k)R^{-1}(k)\mathbf{x}_l(k)e_{j+1}e_{l+1}\right|^q\right\}$$

$$(5.6) \quad \leq C\left[E\left\{\left|\frac{1}{Nk^{1/2}}\sum_{j=K_n}^{n-1}\mathbf{x}'_j(k)R^{-1}(k)\mathbf{x}_j(k)(e_{j+1}^2-\sigma^2)\right|^q\right\}\right.$$

$$+E\left\{\left|\frac{1}{Nk^{1/2}}\sum_{j=K_n}^{n-1}(\mathbf{x}'_j(k)R^{-1}(k)\mathbf{x}_j(k)-k)\right|^q\right\}$$

$$\left.+E\left\{\left|\frac{1}{Nk^{1/2}}\sum_{l=K_n+1}^{n-1}\sum_{j=K_n}^{l-1}\mathbf{x}'_j(k)R^{-1}(k)\mathbf{x}_l(k)e_{j+1}e_{l+1}\right|^q\right\}\right]$$

$$\equiv C\{(\text{I})+(\text{II})+(\text{III})\}.$$

Since $\{\mathbf{x}'_j(k)R^{-1}(k)\mathbf{x}_j(k)(e_{j+1}^2-\sigma^2),\mathcal{M}_{j+1}\}$ is a sequence of martingale differences, where $\mathcal{M}_{j+1}$ is the $\sigma$-algebra generated by $\{e_{j+1},e_j,\ldots\}$, and

$$E\left|\sum_{j=K_n}^{n-1}\mathbf{x}'_j(k)R^{-1}(k)\mathbf{x}_j(k)(e_{j+1}^2-\sigma^2)\right|^q$$

$$\leq E\sup_{K_n\leq m\leq n-1}\left|\sum_{j=K_n}^{m}\mathbf{x}'_j(k)R^{-1}(k)\mathbf{x}_j(k)(e_{j+1}^2-\sigma^2)\right|^q,$$

by Wei ([29], Lemma 2), the assumption that $\sup_{-\infty<t<\infty}E|e_{t+1}^2|^q<\infty$ and the convexity of $x^{q/2}$, $x>0$,

$$(\text{I})\leq C\left(\frac{1}{Nk^{1/2}}\right)^q E\left[\sum_{j=K_n}^{n-1}(\mathbf{x}'_j(k)R^{-1}(k)\mathbf{x}_j(k))^2\right]^{q/2}$$

$$\leq C\left(\frac{1}{Nk}\right)^{q/2}\frac{1}{N}\sum_{j=K_n}^{n-1}E|\mathbf{x}'_j(k)R^{-1}(k)\mathbf{x}_j(k)|^q.$$

Simple algebraic manipulations and Remark 1 yield

$$E|\mathbf{x}'_j(k)R^{-1}(k)\mathbf{x}_j(k)|^q\leq E\|\mathbf{x}_j(k)\|^{2q}\|R^{-1}(k)\|^q\leq CE\|\mathbf{x}_j(k)\|^{2q}$$

$$=CE\left(\sum_{l=1}^k x_{j-l+1}^2\right)^q$$

$$\leq Ck^q k^{-1}\sum_{l=1}^k E|x_{j-l+1}|^{2q}.$$

Since $x_{j-l+1}=\sum_{l_1=0}^\infty b_{l_1}e_{j-l+1-l_1}$, by Wei ([29], Lemma 2), one has, for all integers $j$ and $l$,

$$(5.7) \quad E|x_{j-l+1}|^{2q}\leq C\left(\sum_{l_1=0}^\infty b_{l_1}^2\right)^q\leq C,$$

which further implies, for all $K_n \leq j \leq n$ and all $1 \leq k \leq K_n$,

(5.8) $$E|\mathbf{x}'_j(k)R^{-1}(k)\mathbf{x}_j(k)|^q \leq Ck^q.$$

As a result,

(5.9) $$(\mathrm{I}) \leq C\left(\frac{k}{N}\right)^{q/2}$$

holds for all $1 \leq k \leq K_n$.

Observe that, for all $1 \leq k \leq K_n$,

$$E\left|\frac{1}{N}\sum_{j=K_n}^{n-1}(\mathbf{x}'_j(k)R^{-1}(k)\mathbf{x}_j(k) - k)\right|^q = E|\mathrm{tr}\{R^{-1}(k)(\hat{R}_n(k) - R(k))\}|^q$$

$$= E\left|\sum_{i=1}^{k}\sum_{j=1}^{k} R^{-1}_{i,j}(k)(\hat{r}^{(n)}_{j,i} - r_{j,i})\right|^q$$

$$\leq \left\{\sum_{i=1}^{k}\sum_{j=1}^{k} |R^{-1}_{i,j}(k)|(E|\hat{r}^{(n)}_{j,i} - r_{j,i}|^q)^{1/q}\right\}^q$$

$$\leq C\frac{k^{3q/2}}{N^{q/2}},$$

where $R^{-1}_{i,j}(k), \hat{r}^{(n)}_{i,j}$ and $r_{i,j}$ denote the $(i,j)$ components of $R^{-1}(k), \hat{R}_n(k)$ and $R(k)$, respectively, the first inequality follows from Minkowski's inequality, and the second inequality follows from an inequality given after (2.19) of [17] [which shows that, for all $1 \leq i,j \leq k \leq K_n$, $E|\hat{r}^{(n)}_{j,i} - r_{j,i}|^q \leq CN^{-q/2}$] and the fact that $\sum_{i=1}^{k}\sum_{j=1}^{k}|R^{-1}_{i,j}(k)| \leq Ck^{3/2}$ (see also [21], page 98). Consequently, we have, for all $1 \leq k \leq K_n$,

(5.10) $$(\mathrm{II}) \leq C\frac{k^q}{N^{q/2}}.$$

By Wei ([29], Lemma 2), the moment assumption on $\{e_t\}$ and some algebraic manipulations,

$$(\mathrm{III}) \leq C\left(\frac{1}{Nk^{1/2}}\right)^q E\left[\sum_{l=K_n+1}^{n-1}\left(\sum_{j=K_n}^{l-1}\mathbf{x}'_j(k)R^{-1}(k)\mathbf{x}_l(k)e_{j+1}\right)^2\right]^{q/2}$$

$$\leq C\left\{\left(\frac{1}{Nk^{1/2}}\right)^q E\left[\sum_{l=K_n+1}^{3K_n}\left(\sum_{j=K_n}^{l-1}\mathbf{x}'_j(k)R^{-1}(k)\mathbf{x}_l(k)e_{j+1}\right)^2\right]^{q/2}\right.$$

$$\left. + \left(\frac{1}{Nk^{1/2}}\right)^q\right.$$

$$
(5.11) \quad \times E\left[\sum_{l=3K_n+1}^{n-1}\left(\sum_{j=K_n}^{l-2K_n-1}\mathbf{x}'_j(k)R^{-1}(k)\mathbf{x}_l(k)e_{j+1}\right)^2\right]^{q/2}
$$

$$
+\left(\frac{1}{Nk^{1/2}}\right)^q
$$

$$
\times E\left[\sum_{l=3K_n+1}^{n-1}\left(\sum_{j=l-2K_n}^{l-1}\mathbf{x}'_j(k)R^{-1}(k)\mathbf{x}_l(k)e_{j+1}\right)^2\right]^{q/2}\Bigg\}
$$

$$
\equiv C\{(\text{IV})+(\text{V})+(\text{VI})\}.
$$

By the Cauchy–Schwarz inequality, (5.8), Remark 1 and (5.2), we have, for all $K_n+1 \le l \le 3K_n$,

$$
E\left|\sum_{j=K_n}^{l-1}\mathbf{x}'_j(k)R^{-1}(k)\mathbf{x}_l(k)e_{j+1}\right|^q
$$

$$
(5.12) \quad \le \left\{E|\mathbf{x}'_l(k)R^{-1}(k)\mathbf{x}_l(k)|^q E\left\|\sum_{j=K_n}^{l-1}\mathbf{x}_j(k)e_{j+1}\right\|^{2q}\|R^{-1}(k)\|^q\right\}^{1/2}
$$

$$
\le Ck^q K_n^{q/2},
$$

and for all $3K_n+1 \le l \le n-1$,

$$
(5.13) \quad E\left|\sum_{j=l-2K_n}^{l-1}\mathbf{x}'_j(k)R^{-1}(k)\mathbf{x}_l(k)e_{j+1}\right|^q \le Ck^q K_n^{q/2}.
$$

Hence, by the convexity of $x^{q/2}, x > 0$, (5.12) and (5.13), one obtains, for all $1 \le k \le K_n$,

$$
(5.14) \quad (\text{IV}) \le C\left(\frac{K_n^q k^{q/2}}{N^q}\right) \quad \text{and} \quad (\text{VI}) \le C\left(\frac{K_n k}{N}\right)^{q/2}.
$$

In view of (5.9), (5.10), (5.11) and (5.14), this proof is complete if we can show that (V) is bounded. Observe that

$$
E\left\{\left|\sum_{j=K_n}^{l-2K_n-1}\mathbf{x}'_j(k)R^{-1}(k)\mathbf{x}_l(k)e_{j+1}\right|^q\right\}
$$

$$
\le C\left[E\left\{\left|\sum_{j=K_n}^{l-2K_n-1}\mathbf{x}'_j(k)R^{-1}(k)(\mathbf{x}_l(k)-\tilde{\mathbf{x}}_l^{(n)}(k))e_{j+1}\right|^q\right\}\right.
$$

$$
(5.15)
$$

$$
\left. + E\left\{\left|\sum_{j=K_n}^{l-2K_n-1}\mathbf{x}'_j(k)R^{-1}(k)\tilde{\mathbf{x}}_l^{(n)}(k)e_{j+1}\right|^q\right\}\right]
$$

$$
\equiv C\{(\text{VII})+(\text{VIII})\},
$$

PREDICTION IN AUTOREGRESSIVE PROCESSES 25

where

$$\tilde{\mathbf{x}}_l^{(n)}(k) = \left( \sum_{s=0}^{K_n} b_s e_{l-s}, \ldots, \sum_{s=0}^{K_n} b_s e_{l+1-k-s} \right)',$$

and $b_j$ is defined from the infinite MA representation in Remark 1. Since

$$\mathbf{x}_l(k) - \tilde{\mathbf{x}}_l^{(n)}(k) = \left( \sum_{s=K_n+1}^{\infty} b_s e_{l-s}, \ldots, \sum_{s=K_n+1}^{\infty} b_s e_{l+1-k-s} \right)',$$

reasoning as for (5.8) and (5.12), we have, for all $3K_n + 1 \leq l \leq n-1$, that

(5.16) $$(\text{VII}) \leq C \left( k \sum_{j=K_n+1}^{\infty} b_j^2 \right)^{q/2} N^{q/2} k^{q/2}.$$

It also can be shown that, for all $3K_n + 1 \leq l \leq n-1$,

(5.17)
$$(\text{VIII}) \leq C E \left\{ \left| \left( \sum_{j=K_n}^{l-2K_n-1} \mathbf{x}_j(k) e_{j+1} \right)' \right.\right.$$
$$\left.\left. \times R^{-1}(k) \tilde{R}^{(n)}(k) R^{-1}(k) \left( \sum_{j=K_n}^{l-2K_n-1} \mathbf{x}_j(k) e_{j+1} \right) \right|^{q/2} \right\}$$
$$\leq C \|\tilde{R}^{(n)}(k)\|^{q/2} \|R^{-1}(k)\|^q E \left\| \sum_{j=K_n}^{l-2K_n-1} \mathbf{x}_j(k) e_{j+1} \right\|^q$$
$$\leq C(N^{q/2} k^{q/2}),$$

where $\tilde{R}^{(n)}(k) = E(\tilde{\mathbf{x}}_l^{(n)}(k) \tilde{\mathbf{x}}_l^{(n)'}(k))$, the first inequality is guaranteed by the independence between $\tilde{\mathbf{x}}_l^{(n)}(k)$ and $\sum_{j=K_n}^{l-2K_n-1} \mathbf{x}_j(k) e_{j+1}$, Wei ([29], Lemma 2) and an argument similar to that used for verifying (3.25) and (3.26) of [17]; the third inequality follows from (5.2), Remark 1 and the fact that $\max_{1 \leq k \leq K_n} \tilde{R}^{(n)}(k) < C$, which is ensured by $\sum_{j=0}^{\infty} |b_j| < \infty$. [Note that (5.17) is valid even under a weaker moment assumption, $\sup_{-\infty < t < \infty} E|e_t|^q < \infty, q \geq 2$.] By the convexity of $x^{q/2}, x > 0$ and (5.15)–(5.17), we have, for all $1 \leq k \leq K_n$,

(5.18) $$(\text{V}) \leq C \left( k \sum_{j=K_n+1}^{\infty} b_j^2 + 1 \right)^q.$$

Consequently, the desired property follows from (5.18) and the fact that $n \sum_{j=n+1}^{\infty} b_j^2$ is uniformly bounded. $\square$

REMARK 4. Note that when $q = 2$, Lee and Karagrigoriou [21] introduced a decomposition for

$$E\left(\left\|\frac{1}{\sqrt{Nk}}\sum_{j=K_n}^{n-1}\mathbf{x}_j(k)e_{j+1}\right\|_{R^{-1}(k)}^2 - \sigma^2\right)^q$$

which is similar to those given by (5.6), (5.11) and (5.15). By applying this decomposition, they established (5.5) for the special case of $q = 2$ under (K.1)(c) with $\{e_t\}$ being a sequence of i.i.d. random variables, $E|e_1|^4 < \infty$ and $K_n = o(n^{1/2})$. See Lemma 2.3 of [21] for more details.

As a direct application of Lemma 3, we get the following result.

LEMMA 4. *Assume* (K.1)(b), (K.5), $\sup_{-\infty < t < \infty} E|e_t|^{2q} < \infty$ *for some* $q > 2$ *and* $K_n = O(n^{1/2})$. *Then*

$$(5.19) \quad \lim_{n\to\infty} E\left(\max_{1\le k\le K_n}\left|\left(\left\|\frac{1}{N}\sum_{j=K_n}^{n-1}\mathbf{x}_j(k)e_{j+1}\right\|_{R^{-1}(k)}^2 - \frac{k}{N}\sigma^2\right)L_n^{-1}(k)\right|^q\right) = 0.$$

PROOF. By Lemma 3 and an argument similar to that used for obtaining (3.5) of [27], (5.19) follows. □

Under (K.1)(a) with i.i.d. Gaussian noise, (K.5) and $K_n = o(n^{1/2})$, Shibata ([27], (3.5)) obtained, for any $\epsilon > 0$,

$$(5.20) \quad \lim_{n\to\infty} P\left(\max_{1\le k\le K_n}\left|\frac{\|\hat{\mathbf{a}}_n(k) - \mathbf{a}\|_R^2}{L_n(k)} - 1\right| > \epsilon\right) = 0,$$

where $\hat{\mathbf{a}}_n(k)$ in (5.20) is viewed as an infinite-dimensional vector $(\hat{a}_{1,n}(k), \hat{a}_{2,n}(k), \ldots)'$ with $\hat{a}_{i,n}(k) = 0$ for $i > k$, and $P$ denotes the probability measure. This leads to a lower bound theorem for model selection in independent-realization settings (see Theorem 3.2 of [27]), which serves as an important vehicle for establishing $S_n(k)$'s (and its variants') asymptotic optimality in the sense of (4.11). Recently several authors have attempted to establish (5.20) in non-Gaussian settings. Among them, Lee and Karagrigoriou [21] attained this goal by imposing (K.5), the assumptions described in Remark 4, and for all integer $t$ and all positive integers $k, j$ with $k \ge j$,

$$(5.21) \quad Ew_{t,j}^4(k) \le \bar{C}_1,$$

where $\bar{C}_1 > 0$ is independent of $t, j$ and $k$, and $(w_{t,1}(k), w_{t,2}(k), \ldots, w_{t,k}(k))' = \mathbf{W}_t(k) = R^{-1/2}(k)\mathbf{x}_t(k)$. Their moment assumption on $\{e_t\}$, $E|e_1|^4 < \infty$ (see Remark 4), is considerably weaker than others proposed in the literature. But, (5.21) does not seem to be needed. This is because $w_{t,j}(k)$ is a linear combination in $\{e_l, l \le t\}$, and by an argument similar to that used for showing (5.7), (5.21) holds automatically when the other assumptions are imposed. In fact, by (5.6), (5.9)–(5.11), (5.14)–(5.17) and an idea of Lee and Karagrigoriou ([21], Lemma 2.4), (5.20) can be ensured by a set of weaker assumptions, (K.1)(b), (K.5), $\sup_{-\infty < t < \infty} E|e_t|^4 < \infty$ and $K_n = o(n^{1/2})$. However, to obtain asymptotic expressions for unconditional MSPEs of the least squares predictors with estimated orders, we require a strengthened version of (5.20). The following lemma is given to fulfill this aim.

LEMMA 5. *Let the assumptions of Proposition 2 and* (K.5) *hold. Then, for any* $q > 0$,

$$(5.22) \qquad \lim_{n \to \infty} E\left( \max_{1 \leq k \leq K_n} \left| \frac{\|\hat{\mathbf{a}}_n(k) - \mathbf{a}\|_R^2}{L_n(k)} - 1 \right|^q \right) = 0,$$

*where* $L_n(k)$ *is defined in Proposition 2 of Section 2.*

PROOF. Note that

$$(5.23) \quad \begin{aligned} & \left| \frac{\|\hat{\mathbf{a}}_n(k) - \mathbf{a}\|_R^2}{L_n(k)} - 1 \right| \\ &= \left| \frac{\|\hat{\mathbf{a}}_n(k) - \mathbf{a}(k)\|_{R(k)}^2 - k\sigma^2/N}{L_n(k)} \right| \\ &= \left| \left\| \hat{R}_n^{-1}(k) \frac{1}{N} \sum_{j=K_n}^{n-1} \mathbf{x}_j(k) e_{j+1,k} \right\|_{R(k)}^2 - \frac{k\sigma^2}{N} \right| (L_n(k))^{-1} \\ &= \left| \frac{\mathbf{A}(k) + \mathbf{B}(k) + \mathbf{C}(k) + \mathbf{D}(k)}{L_n(k)} \right|, \end{aligned}$$

where

$$\mathbf{A}(k) = \left\| \frac{1}{N} \sum_{j=K_n}^{n-1} \mathbf{x}_j(k) e_{j+1,k} \right\|_{\hat{R}_n^{-1}(k) R(k) \hat{R}_n^{-1}(k) - R^{-1}(k)}^2,$$

$$\mathbf{B}(k) = \left( \frac{1}{N} \sum_{j=K_n}^{n-1} \mathbf{x}_j'(k) e_{j+1,k} \right) R^{-1}(k) \left( \frac{1}{N} \sum_{j=K_n}^{n-1} \mathbf{x}_j(k) (e_{j+1,k} - e_{j+1}) \right),$$

$$\mathbf{C}(k) = \left( \frac{1}{N} \sum_{j=K_n}^{n-1} \mathbf{x}_j'(k) (e_{j+1,k} - e_{j+1}) \right) R^{-1}(k) \left( \frac{1}{N} \sum_{j=K_n}^{n-1} \mathbf{x}_j(k) e_{j+1} \right)$$

and

$$\mathbf{D}(k) = \left\| \frac{1}{N} \sum_{j=K_n}^{n-1} \mathbf{x}_j(k) e_{j+1} \right\|_{R^{-1}(k)}^2 - \frac{k\sigma^2}{N}.$$

By Lemmas 1 and 2 and Proposition 1, one has, for any $q > 0$ and all $1 \leq k \leq K_n$,

$$E|\mathbf{A}(k)|^q \leq E\left\{\left\|\frac{1}{N}\sum_{j=K_n}^{n-1}\mathbf{x}_j(k)e_{j+1,k}\right\|^{2q}\|\hat{R}_n^{-1}(k)R(k)\hat{R}_n^{-1}(k) - R^{-1}(k)\|^q\right\}$$

(5.24)
$$\leq C\frac{k^{2q}}{N^{3q/2}},$$

$$E|\mathbf{B}(k)|^q \leq E\left\{\left\|\frac{1}{N}\sum_{j=K_n}^{n-1}\mathbf{x}_j(k)e_{j+1,k}\right\|^q\right.$$

(5.25)
$$\left.\times \left\|\frac{1}{N}\sum_{j=K_n}^{n-1}\mathbf{x}_j(k)(e_{j+1,k} - e_{j+1})\right\|^q \|R^{-1}(k)\|^q\right\}$$

$$\leq C\frac{k^q\|\mathbf{a} - \mathbf{a}(k)\|_R^q}{N^q}$$

and

$$E|\mathbf{C}(k)|^q \leq E\left\{\left\|\frac{1}{N}\sum_{j=K_n}^{n-1}\mathbf{x}_j(k)e_{j+1}\right\|^q\right.$$

(5.26)
$$\left.\times \left\|\frac{1}{N}\sum_{j=K_n}^{n-1}\mathbf{x}_j(k)(e_{j+1,k} - e_{j+1})\right\|^q \|R^{-1}(k)\|^q\right\}$$

$$\leq C\frac{k^q\|\mathbf{a} - \mathbf{a}(k)\|_R^q}{N^q}.$$

Now, according to (5.23)–(5.26) and Lemma 3, we have, for sufficiently large $q$,

$$\sum_{k=1}^{K_n} E\left|\frac{\|\hat{\mathbf{a}}_n(k) - \mathbf{a}\|_R^2}{L_n(k)} - 1\right|^q$$

(5.27)
$$\leq C\sum_{k=1}^{K_n}\left(\frac{k^{2q}}{N^{3q/2}L_n^q(k)} + \frac{k^q\|\mathbf{a} - \mathbf{a}(k)\|_R^q}{N^q L_n^q(k)} + \frac{k^{q/2}}{N^q L_n^q(k)}\right)$$

$$\leq C\sum_{k=1}^{K_n}\left(\frac{k^q}{N^{q/2}} + \sum_{k=1}^{k_n^*}\frac{k^{q/2}}{k_n^{*q}} + \sum_{k=k_n^*+1}^{K_n}k^{-q/2}\right) = o(1),$$

where the second inequality is ensured by

$$\frac{\|\mathbf{a} - \mathbf{a}(k)\|_R^2}{L_n(k)} \leq 1$$

and $NL_n(k) \geq C\max\{k, k_n^*\}$, and the equality is guaranteed by $k_n^* \to \infty$. Consequently, (5.22) is ensured by (5.27) and Jensen's inequality. □

Under (K.1)(a) with i.i.d. Gaussian noise, Shibata ([27], Lemmas 4.2, 4.3 and 4.4 and Proposition 4.1) established that, for all $1 \leq k, j \leq K_n$ with $K_n \leq n - 1$,

$$(5.28) \quad E|S_k^2 - \sigma_k^2 - (S_j^2 - \sigma_j^2)|^4 \leq C(\max\{k,j\})^{1/2} N^{-2} \|\mathbf{a}(j) - \mathbf{a}(k)\|^4,$$

where $S_k^2 = (1/N) \sum_{t=K_n}^{n-1} e_{t+1,k}^2$, $\mathbf{a}(j)$ and $\mathbf{a}(k)$ are viewed as infinite-dimensional vectors and $\sigma_k^2 = E(S_k^2)$. [Since, by Remark 1, $D_l \|\mathbf{a}(j) - \mathbf{a}(k)\|^2 \leq \|\mathbf{a}(j) - \mathbf{a}(k)\|_R^2 \leq D_u \|\mathbf{a}(j) - \mathbf{a}(k)\|^2$ for some $0 < D_l \leq D_u < \infty$ independent of $j$ and $k$, $\|\mathbf{a}(j) - \mathbf{a}(k)\|^4$ in (5.28) can be replaced by $\|\mathbf{a}(j) - \mathbf{a}(k)\|_R^4$.] However, his approach, based on heavy calculations for the cross moments of Gaussian random variables, is not easy to extend to higher-moment cases. Moreover, the term $(\max\{k,j\})^{1/2}$ is cumbersome because it is difficult to infer how this term varies with the exponent on the left-hand side of (5.28). Lemma 6 below not only shows that this term is not needed, but also provides a general bound valid for each $q \geq 2$. For some applications of Lemma 6, see Remark 6 and the proofs of Theorems 2 and 3.

LEMMA 6. *Assume* (K.1)(a), $\sup_{-\infty < t < \infty} E|e_t|^{2q} < \infty$ *with* $q \geq 2$ *and* $K_n \leq n - 1$. *Then*

$$(5.29) \quad \max_{1 \leq k \leq K_n} E|S_k^2 - \sigma_k^2|^q \leq CN^{-q/2},$$

*and for all* $1 \leq k, j \leq K_n$,

$$(5.30) \quad E|S_k^2 - \sigma_k^2 - (S_j^2 - \sigma_j^2)|^q \leq CN^{-q/2}\|\mathbf{a}(j) - \mathbf{a}(k)\|_R^q.$$

PROOF. First observe that

$$|S_k^2 - \sigma_k^2|^q = \left|\frac{1}{N}\sum_{t=K_n}^{n-1} e_{t+1,k}^2 - E(e_{t+1,k}^2)\right|^q.$$

Since, according to Remark 1 and the definition of $e_{t,k}$ [given after (2.2)], $e_{t,k}$, $t = \ldots, -1, 0, 1, \ldots$, is a linear process, by Findley and Wei ([11], the first moment bound theorem) one has, for all $1 \leq k \leq K_n$,

$$(5.31) \quad \begin{aligned} E|S_k^2 - \sigma_k^2|^q &\leq C\frac{1}{N^q}\left\{\sum_{i=K_n}^{n-1}\sum_{j=K_n}^{n-1}[E(e_{i+1,k}e_{j+1,k})]^2\right\}^{q/2} \\ &\leq C\frac{1}{N^{q/2}}\left(\sum_{j=-\infty}^{\infty} \xi_j^{*2}\right)^{q/2}, \end{aligned}$$

where $\xi^*_{i-j} = E(e_{i+1,k} e_{j+1,k})$. By Remark 1,

$$\sum_{j=-\infty}^{\infty} \xi^{*2}_j = 2\pi \int_{-\pi}^{\pi} \left( \left| \sum_{j=0}^{k} a_j(k) e^{-ij\lambda} \right|^2 f(\lambda) \right)^2 d\lambda$$

$$\leq 4\pi^2 f_2^2 \left( \sum_{j=0}^{k} |a_j(k)| \right)^4,$$

where $a_0(k) = 1$. By Berk ([4], Lemma 4),

$$\sup_{0 \leq k < \infty} \sum_{j=0}^{k} |a_j(k)| < \infty.$$

This fact and (5.31) yield (5.29).

To obtain (5.30), assume $k < j$. Then

(5.32)
$$|S_k^2 - \sigma_k^2 - (S_j^2 - \sigma_j^2)|$$
$$= \left| \frac{1}{N} \sum_{t=K_n}^{n-1} e_{t+1,k}^2 - \|\mathbf{a}(j) - \mathbf{a}(k)\|_R^2 - e_{t+1,j}^2 \right|$$
$$\leq \left| \frac{1}{N} \sum_{t=K_n}^{n-1} (e_{t+1,k} - e_{t+1,j}) e_{t+1,k} - \|\mathbf{a}(j) - \mathbf{a}(k)\|_R^2 \right|$$
$$+ \left| \frac{1}{N} \sum_{t=K_n}^{n-1} (e_{t+1,k} - e_{t+1,j}) e_{t+1,j} \right|.$$

By Findley and Wei ([11], the first moment bound theorem) again, the expectation of the first term on the right-hand side of (5.32) is bounded by

$$CN^{-q} \left( \sum_{t_1=K_n}^{n-1} \sum_{t_2=K_n}^{n-1} E((e_{t_1+1,k} - e_{t_1+1,j})(e_{t_2+1,k} - e_{t_2+1,j})) \xi^*_{t_1-t_2} \right)^{q/2}$$

$$\leq CN^{-q/2} \{E(e_{1,k} - e_{1,j})^2\}^{q/2} \left( \sum_{l=-\infty}^{\infty} |\xi^*_l| \right)^{q/2} \leq CN^{-q/2} \|\mathbf{a}(j) - \mathbf{a}(k)\|_R^q,$$

where the second inequality follows from $E(e_{1,k} - e_{1,j})^2 = \|\mathbf{a}(j) - \mathbf{a}(k)\|_R^2$ and

$$\sum_{j=-\infty}^{\infty} |\xi^*_j| \leq C \left( \sum_{i=0}^{k} |a_i(k)| \sum_{j=0}^{\infty} |b_j| \right)^2 < C.$$

Similarly, an upper bound for the expectation of the second term on the right-hand side of (5.32) is also given by $CN^{-q/2}\|\mathbf{a}(j) - \mathbf{a}(k)\|_R^q$. As a result, (5.30) holds. □

REMARK 5. Assume $\{e_t\}$ in Lemma 6 is a sequence of martingale differences corresponding to an increasing sequence of $\sigma$-fields of events, $\{\mathcal{F}_t\}$. Further, assume that

$$E(e_t^2|\mathcal{F}_{t-1}) = \sigma^2 \quad \text{a.s.}$$

for $t = \ldots, -1, 0, 1, \ldots$, and

$$\sup_{-\infty < t < \infty} E(|e_t|^{2q}|\mathcal{F}_{t-1}) \leq C < \infty \quad \text{a.s.}$$

Then by the same argument as in Lemma 6, (5.29) and (5.30) are still valid.

According to the decomposition of $S_n(k)$ given by (4.1) of [27],

$$\begin{aligned}
S_n(k) = NL_n(k) + 2k(\hat{\sigma}_k^2 - \sigma^2) + (k\sigma^2 - N\|\hat{a}(k) - a(k)\|_{\hat{R}_n(k)}^2) \\
+ N\sigma^2 + N(S_k^2 - \sigma_k^2).
\end{aligned} \tag{5.33}$$

Equality (5.33) yields

$$P(\hat{k}_n^S = k) \leq P(S_n(k) \leq S_n(k_n^*)) = P\left(\frac{S_n(k)}{NL_n(k)} \leq \frac{S_n(k_n^*)}{NL_n(k)}\right)$$

$$\leq \sum_{i=1}^{5} P(|V_{in}(k)| \geq (1/5)V_n(k)), \tag{5.34}$$

where

$$V_{1n}(k) = -\frac{2k(\hat{\sigma}_k^2 - \sigma^2)}{NL_n(k)},$$

$$V_{2n}(k) = \frac{2k_n^*(\hat{\sigma}_{k_n^*}^2 - \sigma^2)}{NL_n(k)},$$

$$V_{3n}(k) = -\frac{k\sigma^2 - N\|\hat{\mathbf{a}}_n(k) - \mathbf{a}(k)\|_{\hat{R}_n(k)}^2}{NL_n(k)},$$

$$V_{4n}(k) = \frac{k_n^*\sigma^2 - N\|\hat{\mathbf{a}}_n(k_n^*) - \mathbf{a}(k_n^*)\|_{\hat{R}_n(k_n^*)}^2}{NL_n(k)},$$

$$V_{5n}(k) = -\frac{S_k^2 - \sigma_k^2 - S_{k_n^*}^2 - \sigma_{k_n^*}^2}{L_n(k)}$$

and

$$V_n(k) = \frac{L_n(k) - L_n(k_n^*)}{L_n(k)}.$$

By (5.34), Chebyshev's inequality and Lemmas 1–6, we can obtain an upper bound for $P(\hat{k}_n^S = k)$ through moment bounds for $V_{in}(k)/V_n(k), i = 1, \ldots, 5$. This is our motivation for verifying that, for any $r > 0$,

(5.35) $$\lim_{n \to \infty} E\left(\frac{L_n(\hat{k}_n^S)}{L_n(k_n^*)} - 1\right)^r = 0;$$

see Theorem 3 below for more details. Equality (5.35), which provides a (general) moment convergence result for $L_n(\hat{k}_n^S)/L_n(k_n^*)$, is the key to proving Theorem 1, and can be applied to verify (5.74) and (5.75), which are the final steps in the proof of Theorem 2. As a byproduct, (5.35) also yields $\hat{k}_n^S$'s asymptotic efficiency for independent-realization predictions in the sense of Shibata [27], as defined in (4.11). To see this, first notice that, for any random variable $\hat{I}_n \in \mathcal{J}_n$ [$\mathcal{J}_n$ is defined after (4.6)],

$$E\{(y_{n+1} - \hat{y}_{n+1}(\hat{I}_n))^2 | x_1, \ldots, x_n\} - \sigma^2 = \|\hat{\mathbf{a}}_n(\hat{I}_n) - \mathbf{a}\|_R^2.$$

It is also not difficult to show that

(5.36) $$1 - \max_{1 \leq k \leq K_n} \left|\frac{\|\hat{\mathbf{a}}_n(k) - \mathbf{a}\|_R^2}{L_n(k)} - 1\right| \leq \inf_{\hat{I}_n \in \mathcal{J}_n} \frac{\|\hat{\mathbf{a}}_n(\hat{I}_n) - \mathbf{a}\|_R^2}{L_n(\hat{I}_n)}$$
$$\leq \frac{\inf_{\hat{I}_n \in \mathcal{J}_n} \|\hat{\mathbf{a}}_n(\hat{I}_n) - \mathbf{a}\|_R^2}{L_n(k_n^*)}$$
$$\leq \frac{\|\hat{\mathbf{a}}_n(k_n^*) - \mathbf{a}\|_R^2}{L_n(k_n^*)}.$$

Since by (5.20) both sides of (5.36) converge to 1 in probability, (4.11) can be rewritten as

(5.37) $$\frac{L_n(\hat{k}_n)}{L_n(k_n^*)} - 1 = o_p(1) \qquad \left(\text{or } \frac{\|\hat{\mathbf{a}}_n(\hat{k}_n) - \mathbf{a}\|_R^2}{L_n(k_n^*)} - 1 = o_p(1)\right).$$

Obviously,

(5.38) $$\frac{L_n(\hat{k}_n^S)}{L_n(k_n^*)} - 1 = o_p(1)$$

is an immediate consequence of (5.35).

THEOREM 3. *Let the assumptions of Proposition 2 and* (K.5) *hold. Then* (5.35) *holds.*

PROOF. Let $\epsilon > 0$. By (5.34),

$$E\left(\frac{L_n(\hat{k}_n^S)}{L_n(k_n^*)} - 1\right)^r$$

$$= \sum_{k=1}^{K_n} \left(\frac{L_n(k)}{L_n(k_n^*)} - 1\right)^r P(\hat{k}_n^S = k)$$

(5.39)

$$\leq \epsilon^r + \sum_{k \in A_{\epsilon,n}} \left(\frac{L_n(k)}{L_n(k_n^*)} - 1\right)^r P(\hat{k}_n^S = k)$$

$$\leq \epsilon^r + \sum_{i=1}^{5} \left\{\sum_{k \in A_{\epsilon,n}} \left(\frac{L_n(k)}{L_n(k_n^*)} - 1\right)^r P(|V_{in}(k)| > (1/5)V_n(k))\right\},$$

where

(5.40) $$A_{\epsilon,n} = \left\{k : 1 \leq k \leq K_n, \frac{L_n(k)}{L_n(k_n^*)} - 1 > \epsilon\right\}.$$

To obtain (5.35), it suffices to show that the value inside the braces of (5.39) is asymptotically negligible in the sense that

(5.41) $$\lim_{n \to \infty} \sum_{k \in A_{\epsilon,n}} \left(\frac{L_n(k)}{L_n(k_n^*)} - 1\right)^r P(|V_{in}(k)| > (1/5)V_n(k)) = 0,$$

for $i = 1, \ldots, 5$.

By (4.2) of [27] and some algebraic manipulations,

(5.42)
$$|\hat{\sigma}_k^2 - \sigma^2| \leq \|\hat{\mathbf{a}}_n(k) - \mathbf{a}(k)\|_{\hat{R}_n(k)}^2 + \|\mathbf{a}(k) - \mathbf{a}\|_R^2 + |S_k^2 - \sigma_k^2|$$

$$\leq (\|\hat{R}_n(k) - R(k)\|\|R^{-1}(k)\| + 1)\|\hat{\mathbf{a}}_n(k) - \mathbf{a}\|_R^2 + |S_k^2 - \sigma_k^2|.$$

By Ing and Wei ([17], Lemma 2), (5.22), (5.29), (5.42) and Hölder's inequality, one has, for any $q > 0$ and all $1 \leq k \leq K_n$,

(5.43) $$E|V_{1n}(k)|^q \leq C\left(\frac{k^q}{N^q} + N^{-q/2}\right).$$

Since

$$\frac{L_n(k)}{L_n(k_n^*)} \leq \begin{cases} CL_n^{-1}(k_n^*), & \text{if } 1 \leq k \leq k_n^*, \\ Ck, & \text{if } k_n^* + 1 \leq k \leq K_n \end{cases}$$

and

$$V_n^{-1}(k) \leq \frac{1+\epsilon}{\epsilon} \qquad \text{for } k \in A_{\epsilon,n},$$

by the Chebyshev inequality and (5.43), we have, for large $q$,

$$\sum_{k\in A_{\epsilon,n}} \left(\frac{L_n(k)}{L_n(k_n^*)} - 1\right)^r P(|V_{1n}(k)| > (1/5)V_n(k))$$

$$\leq C \sum_{k\in A_{\epsilon,n}} \left(\frac{L_n(k)}{L_n(k_n^*)}\right)^r V_n^{-(q-r)}(k)\left(\frac{k^q}{N^q} + N^{-q/2}\right)$$

(5.44)
$$\leq C\left(\frac{1+\epsilon}{\epsilon}\right)^{q-r} \left\{\sum_{k=1}^{k_n^*}\left(\frac{k^q}{k_n^{*r}N^{q-r}} + \frac{1}{k_n^{*r}N^{q/2-r}}\right)\right.$$

$$\left. + \sum_{k=k_n^*+1}^{K_n}\left(\frac{k^{q+r}}{N^q} + \frac{k^r}{N^{q/2}}\right)\right\}$$

$$= o(1).$$

Therefore, (5.41) holds for $i = 1$. Similarly, we also have, for any $q > 0$ and all $1 \leq k \leq K_n$,

(5.45)
$$E|V_{2n}(k)|^q \leq C\left\{\frac{k_n^{*q}}{N^q} + N^{-q/2}\right\}.$$

By (5.45) and the same argument used for verifying (5.44), (5.41) holds for $i = 2$.

For $i = 3$, we have

(5.46)
$$|V_{3n}(k)| \leq \left|1 - \frac{\|\hat{\mathbf{a}}_n(k) - \mathbf{a}\|_R^2}{L_n(k)}\right|$$

$$+ \left|\frac{\|\hat{\mathbf{a}}_n(k) - \mathbf{a}(k)\|_{\hat{R}_n(k)}^2 - \|\hat{\mathbf{a}}_n(k) - \mathbf{a}(k)\|_{R(k)}^2}{L_n(k)}\right|$$

$$\leq \left|1 - \frac{\|\hat{\mathbf{a}}_n(k) - \mathbf{a}\|_R^2}{L_n(k)}\right|$$

$$+ \frac{\|\hat{R}_n(k) - R(k)\|\|R^{-1}(k)\|\|\hat{\mathbf{a}}_n(k) - \mathbf{a}\|_R^2}{L_n(k)}.$$

By Ing and Wei ([17], Lemma 2) and (5.23)–(5.26), we have, for any $q > 0$ and all $1 \leq k \leq K_n$,

(5.47)
$$E|V_{3n}(k)|^q \leq C\left(\frac{k^q}{N^{q/2}} + \frac{k^{q/2}}{N^q L_n^q(k)}\right).$$

Now, by taking a sufficiently large $q$,

$$\sum_{k\in A_{\epsilon,n}} \left(\frac{L_n(k)}{L_n(k_n^*)} - 1\right)^r P(|V_{3n}(k)| > (1/5)V_n(k))$$

$$\leq C \sum_{k \in A_{\epsilon,n}} \left(\frac{L_n(k)}{L_n(k_n^*)}\right)^r V_n^{-(q-r)}(k) \left(\frac{k^q}{N^{q/2}} + \frac{k^{q/2}}{N^q L_n^q(k)}\right)$$

$$\leq C \left(\frac{1+\epsilon}{\epsilon}\right)^{q-r} \left\{\sum_{k=1}^{k_n^*} \left(\frac{k^q}{k_n^{*r} N^{q/2-r}} + \frac{k^{q/2}}{k_n^{*q}}\right) \right.$$

$$\left. + \sum_{k=k_n^*+1}^{K_n} \left(\frac{k^{q+r}}{N^{q/2}} + k_n^{*-r} k^{-q/2+r}\right)\right\}$$

$$= o(1).$$

Consequently, (5.41) holds for $i=3$. Similarly, we can also show that, for any $q > 0$ and all $1 \leq k \leq K_n$,

(5.48) $$E|V_{4n}(k)|^q \leq C \left(\frac{k_n^{*q}}{N^{q/2}} + \frac{k_n^{*q/2}}{N^q L_n^q(k)}\right),$$

and, hence, (5.41) holds for $i=4$.

Since (K.3) is assumed, by (5.30) one has, for any $q > 0$ and all $1 \leq k \leq K_n$,

(5.49) $$E|V_{5n}(k)|^q \leq C \frac{\|\mathbf{a}(k) - \mathbf{a}(k_n^*)\|_R^q}{N^{q/2} L_n^q(k)}.$$

This gives for large $q$,

$$\sum_{k \in A_{\epsilon,n}} \left(\frac{L_n(k)}{L_n(k_n^*)} - 1\right)^r P(|V_{5n}(k)| > (1/5) V_n(k))$$

$$\leq C \left(\frac{1+\epsilon}{\epsilon}\right)^{q-r} \sum_{k=1}^{K_n} \frac{\|\mathbf{a}(k) - \mathbf{a}(k_n^*)\|_R^q}{L_n^r(k_n^*) N^{q/2} L_n^{q-r}(k)}$$

$$\leq C \left(\frac{1+\epsilon}{\epsilon}\right)^{q-r}$$

$$\times \left(\sum_{k=1}^{k_n^*} \frac{\|\mathbf{a} - \mathbf{a}(k)\|_R^q}{L_n^r(k_n^*) N^{q/2} L_n^{q-r}(k)} + \sum_{k=k_n^*+1}^{K_n} \frac{\|\mathbf{a} - \mathbf{a}(k_n^*)\|_R^q}{L_n^{q/2}(k_n^*) N^{q/2} L_n^{q/2}(k)}\right)$$

$$\leq C \left(\frac{1+\epsilon}{\epsilon}\right)^{q-r} \left(\sum_{k=1}^{k_n^*} k_n^{*-q/2} + \sum_{k=k_n^*+1}^{K_n} k^{-q/2}\right)$$

$$= o(1).$$

Hence, (5.41) holds for $i=5$. Consequently (5.35) follows. $\square$

REMARK 6. In this remark we show that (5.38) can be directly verified [without the help of (5.35)] under much weaker conditions than those of Theorem 3. [Recall that the main purpose of Theorem 3 is to provide a general moment bound for $L_n(\hat{k}_n^S)/L_n(k_n^*)$, and (5.38) is only its byproduct.] To see this, first notice that, by assuming (K.1)(a), $k_n^* \to \infty$ as $n \to \infty$ and the Gaussianity of $\{e_t\}$, Shibata ([27], Proposition 4.1) proved that

$$\max_{1 \leq k \leq K_n} |V_{5n}(k)| = o_p(1). \tag{5.50}$$

As can be seen from [27] and [21], (5.50) and (5.20) are the two most important tools for verifying (5.38). However, by (5.30), (5.49) and the assumption that $k_n^* \to \infty$ as $n \to \infty$, one has, for $q > 2$,

$$\sum_{k=1}^{K_n} E|V_{5n}(k)|^q \leq \sum_{k=1}^{k_n^*} \{NL_n(k_n^*)\}^{-q/2} + \sum_{k=k_n^*+1}^{K_n} \{NL_n(k)\}^{-q/2} = o(1),$$

which in turns implies (5.50). As a result, (5.50) still follows if the Gaussian assumption on $\{e_t\}$ is replaced by $\sup_{-\infty < t < \infty} E|e_t|^q < \infty$, $q > 4$. This fact, a result given before Lemma 5 [which shows that (5.20) can be guaranteed by a set of rather mild assumptions], (5.29) and the same argument as the one given in Theorem 4.1 of [27] together yield that (K.1)(b), (K.5), $K_n = o(n^{1/2})$ and $\sup_{-\infty < t < \infty} E|e_t|^q < \infty$, $q > 4$, are sufficient to confirm (5.38). Recently, under model (1.1) with i.i.d. but non-Gaussian noise, Bhansali [5], Karagrigoriou [19] and Lee and Karagrigoriou [21] also obtained (5.38). However, all these papers required a more stringent moment assumption, $E|e_1|^q < \infty$ with $q$ larger than or equal to 8; see Section 6 for more discussion.

The following corollary deals with the moment properties of $L_n(\hat{k}_n)$ with $\hat{k}_n = \hat{k}_n^A, \hat{k}_n^F, \hat{k}_n^{S_p}$ and $\hat{k}_n^C$.

COROLLARY 2. *Let the assumptions of Theorem 3 hold. Then* (5.35) *holds with $\hat{k}_n^S$ replaced by $\hat{k}_n^A, \hat{k}_n^F, \hat{k}_n^{S_p}$ or $\hat{k}_n^C$.*

PROOF. Define $G_n^{(1)}(k) = S_n(k) - N \exp(\text{AIC}(k))$, $G_n^{(2)}(k) = S_n(k) - N(\text{FPE}(k))$, $G_n^{(3)}(k) = S_n(k) - N(\text{S}_p(k))$ and $G_n^{(4)}(k) = S_n(k) - C_p(k)$. By arguments similar to those in Theorem 3 and Shibata ([27], Theorem 4.2), (5.41) still holds with $|V_{in}(k)|$ replaced by

$$\frac{|G_n^{(j)}(k) - G_n^{(j)}(k_n^*)|}{NL_n(k)},$$

$j \in \{1, 2, 3, 4\}$, or with $1/5$ replaced by any positive number independent of $n$. Viewing the proof of Theorem 3, the claimed properties are guaranteed by this finding. □

We are now in position to prove Theorem 1.

PROOF OF THEOREM 1. First note that

(5.51)
$$\frac{E(y_{n+1} - \hat{y}_{n+1}(\hat{k}_n))^2 - \sigma^2}{L_n(k_n^*)}$$
$$= E\bigg\{\bigg(\frac{\|\hat{\mathbf{a}}_n(\hat{k}_n) - \mathbf{a}\|_R^2}{L_n(\hat{k}_n)} - 1\bigg)\frac{L_n(\hat{k}_n)}{L_n(k_n^*)}\bigg\} + E\bigg(\frac{L_n(\hat{k}_n)}{L_n(k_n^*)}\bigg),$$

where $\hat{k}_n = \hat{k}_n^S, \hat{k}_n^A, \hat{k}_n^F, \hat{k}_n^{S_p}$ or $\hat{k}_n^C$. By (5.22), Theorem 3 and Corollary 2, the first expectation on the right-hand side of (5.51) converges to 0, whereas the second expectation converges to 1. As a result, Theorem 1 follows. □

To obtain Theorem 2, we still need the following two lemmas, Lemmas 7 and 8.

LEMMA 7. *Suppose that the assumptions of Proposition 2 hold. Then, for any $q > 0$,*

(5.52)
$$\lim_{n\to\infty} E\bigg(\max_{1\le k\le K_n}\bigg|\frac{\mathbf{f}(k) - \mathbf{f}_1(k)}{L_n^{1/2}(k)}\bigg|^q\bigg) = 0,$$

*where $\mathbf{f}(k)$ is defined after (2.3) and for $n \ge \sqrt{n} + K_n + 1$ and $\sqrt{n} \ge 2K_n$,*

$$\mathbf{f}_1(k) = \mathbf{x}_n^{*\prime}(k)R^{-1}(k)\frac{1}{N}\sum_{j=K_n}^{n-\sqrt{n}-1}\mathbf{x}_j(k)e_{j+1}$$

*with $\mathbf{x}_n^{*\prime}(k) = (x_n^*, \ldots, x_{n-k+1}^*) = (\sum_{r=0}^{\sqrt{n}/2-K_n} b_r e_{n-r}, \ldots, \sum_{r=0}^{\sqrt{n}/2-K_n} b_r e_{n-k+1-r}).$*

PROOF. By (K.4) we can assume that $n \ge \sqrt{n} + K_n + 1$ and $\sqrt{n} \ge 2K_n$ without loss of generality. To obtain (5.52), we first show that, for any $q > 0$,

(5.53)
$$\lim_{n\to\infty} E\bigg(\max_{1\le k\le K_n}\bigg|\frac{\mathbf{f}(k) - \mathbf{f}_0(k)}{L_n^{1/2}(k)}\bigg|^q\bigg) = 0,$$

where

$$\mathbf{f}_0(k) = \mathbf{x}_n^{*\prime}(k)\hat{R}_n^{\circ-1}(k)\frac{1}{N}\sum_{j=K_n}^{n-\sqrt{n}-1}\mathbf{x}_j(k)e_{j+1,k}$$

with

$$\hat{R}_n^{\circ}(k) = \frac{1}{N}\sum_{j=K_n}^{n-\sqrt{n}-1}\mathbf{x}_j(k)\mathbf{x}_j'(k).$$

Set $q \geq 2/3$. By Proposition 1, Lemmas 1 and 2, Hölder's inequality and an argument similar to that used for obtaining (3.13) of [17],

$$E\left(\max_{1 \leq k \leq K_n} \left|\frac{\mathbf{f}(k) - \mathbf{f}_0(k)}{L_n^{1/2}(k)}\right|^q\right)$$

$$\leq C \sum_{i=1}^{K_n} L_n^{-q/2}(i)$$

$$\times \left\{\left[E(\|\mathbf{x}_n(i) - \mathbf{x}_n^*(i)\|^{3q})E(\|\hat{R}_n^{-1}(i)\|^{3q})\right.\right.$$

$$\left.\times E\left(\left\|N^{-1} \sum_{j=K_n}^{n-1} \mathbf{x}_j(i)e_{j+1,i}\right\|^{3q}\right)\right]^{1/3}$$

(5.54) $\qquad + \left[E(\|\mathbf{x}_n^*(i)\|^{3q})E(\|\hat{R}_n^{-1}(i) - \hat{R}_n^{\circ^{-1}}(i)\|^{3q})\right.$

$$\left.\times E\left(\left\|N^{-1} \sum_{j=K_n}^{n-1} \mathbf{x}_j(i)e_{j+1,i}\right\|^{3q}\right)\right]^{1/3}$$

$$+ \left[E(\|\mathbf{x}_n^*(i)\|^{3q})E(\|\hat{R}_n^{\circ^{-1}}(i)\|^{3q})\right.$$

$$\left.\left.\times E\left(\left\|N^{-1} \sum_{j=n-\sqrt{n}}^{n-1} \mathbf{x}_j(i)e_{j+1,i}\right\|^{3q}\right)\right]^{1/3}\right\}$$

$$\leq C \sum_{i=1}^{K_n} L_n^{-q/2}(i)\left[\left(i \sum_{j \geq \sqrt{n}/2 - K_n + 1}^{\infty} b_j^2\right)^{q/2} i^{q/2} N^{-q/2}\right.$$

$$\left.+ i^{2q}N^{-5q/4} + i^q N^{-3q/4}\right].$$

By observing that $L_n^{-q/2}(i) \leq (i/N)^{-q/2}$, the right-hand side of (5.54) is bounded by

(5.55) $\quad C \sum_{i=1}^{K_n}\left[\left(i \sum_{j \geq \sqrt{n}/2 - K_n + 1}^{\infty} b_j^2\right)^{q/2} + i^{3q/2}N^{-3q/4} + i^{q/2}N^{-q/4}\right].$

Moreover, since (K.1)(b) implies that $n \sum_{j \geq n} b_j^2 = o(1)$, by taking sufficiently large $q$, (5.55) converges to 0, and hence (5.53) holds for sufficiently large $q$. This result and Jensen's inequality further guarantee that (5.53) is valid for any $q > 0$.

In view of (5.53), (5.52) follows if one can show that, for any $q > 0$,

(5.56) $$\lim_{n \to \infty} E\left(\max_{1 \leq k \leq K_n} \left|\frac{\mathbf{f}_0(k) - \mathbf{f}_1(k)}{L_n^{1/2}(k)}\right|^q\right) = 0.$$

By Proposition 1, Lemmas 1 and 2 and an argument similar to that used for showing (3.21) of [17], one has, for sufficiently large $q$,

$$E\left(\max_{1 \leq k \leq K_n} \left|\frac{\mathbf{f}_0(k) - \mathbf{f}_1(k)}{L_n^{1/2}(k)}\right|^q\right) \leq \sum_{i=1}^{K_n} E\left|\frac{\mathbf{f}_0(i) - \mathbf{f}_1(i)}{L_n^{1/2}(i)}\right|^q \leq C \sum_{i=1}^{K_n} \frac{i^q}{N^{q/2}} = o(1),$$

which, together with Jensen's inequality, yields (5.56). □

LEMMA 8. *Assume* (K.1)(a), $\sup_{-\infty < t < \infty} E|e_t|^{2q} < \infty, q \geq 1$, *and* $K_n = o(n^{1/2})$. *Then, for* $n \geq \sqrt{n} + K_n + 1$ *and* $\sqrt{n} \geq 2K_n$,

(5.57) $$\max_{\substack{1 \leq i, l \leq K_n \\ i \neq l}} \frac{E|\mathbf{f}_1(i) - \mathbf{f}_1(l)|^{2q}}{|(i-l)/N|^q} \leq C.$$

PROOF. Without loss of generality, assume that $1 \leq i < l \leq K_n$. First observe that

(5.58) $$\begin{aligned}\mathbf{f}_1(i) - \mathbf{f}_1(l) &= \mathbf{x}_n^{*\prime}(i) R^{-1}(l) \frac{1}{N} \sum_{j=K_n}^{n-\sqrt{n}-1} (\mathbf{x}_j(i) - \mathbf{x}_j(l)) e_{j+1} \\ &\quad + (\mathbf{x}_n^*(i) - \mathbf{x}_n^*(l))' R^{-1}(l) \frac{1}{N} \sum_{j=K_n}^{n-\sqrt{n}-1} \mathbf{x}_j(l) e_{j+1} \\ &\quad + \mathbf{x}_n^{*\prime}(i)(R^{-1}(i) - R^{-1}(l)) \frac{1}{N} \sum_{j=K_n}^{n-\sqrt{n}-1} \mathbf{x}_j(i) e_{j+1} \\ &\equiv (\mathrm{I}) + (\mathrm{II}) + (\mathrm{III}),\end{aligned}$$

where $\mathbf{x}_n^*(i)$, $\mathbf{x}_j(i)$ and $R^{-1}(i)$, respectively, are regarded as $l$-dimensional vectors and an $l \times l$ matrix with undefined entries set to 0. To obtain (5.57), it suffices to show that, for all $1 \leq i < l \leq K_n$,

$$E|(T)|^{2q} \leq C\left(\frac{l-i}{N}\right)^q,$$

for all $T = (\mathrm{I})$, $(\mathrm{II})$ and $(\mathrm{III})$. (Note that as mentioned before Proposition 1, $C$ is used to denote some positive number independent of $i, l, K_n$ and $n$.)

Since $\mathbf{x}_n^*(i)$ is independent of the remaining part of (I), by an argument similar to that used for showing (5.17), one has, for all $1 \leq i < l \leq K_n$,

$$E(|(\mathrm{I})|^{2q}) \leq CE \left\| \frac{1}{N} \sum_{j=K_n}^{n-\sqrt{n}-1} (\mathbf{x}_j(i) - \mathbf{x}_j(l)) e_{j+1} \right\|_{R^{-1}(l) R_*^{(n)}(i) R^{-1}(l)}^{2q}$$

(5.59)
$$\leq C \|R^{-1}(l)\|^{2q} \|R_*^{(n)}(i)\|^q$$

$$\times E \left\| \frac{1}{N} \sum_{j=K_n}^{n-\sqrt{n}-1} (\mathbf{x}_j(i) - \mathbf{x}_j(l)) e_{j+1} \right\|^{2q}$$

$$\leq C \left( \frac{l-i}{N} \right)^q,$$

where $R_*^{(n)}(i) = E(\mathbf{x}_n^*(i) \mathbf{x}_n^{*'}(i))$ and the last inequality follows from Lemma 2 and

$$\max_{1 \leq k \leq K_n} \|R_*^{(n)}(k)\| < C,$$

which is guaranteed by $\sum_{i=0}^{\infty} |b_i| < \infty$. Similarly, for all $1 \leq i < l \leq K_n$,

$$E(|(\mathrm{II})|^{2q}) \leq CE \left\| \frac{1}{N} \sum_{j=K_n}^{n-\sqrt{n}-1} \mathbf{x}_j(l) e_{j+1} \right\|_{R^{-1}(l) D_n(l,i) R^{-1}(l)}^{2q}$$

$$\leq CE \left\| \frac{1}{N} \sum_{j=K_n}^{n-\sqrt{n}-1} \mathbf{z}_j^{(n)}(l,i) e_{j+1} \right\|^{2q},$$

where

$$D_n(l,i) = \begin{pmatrix} 0_{i \times i} & 0_{i \times (l-i)} \\ 0_{(l-i) \times i} & R_*^{(n)}(l-i) \end{pmatrix}$$

and

$$\mathbf{z}_j^{(n)}(l,i) = (z_{j,1}^{(n)}(l,i), \ldots, z_{j,l-i}^{(n)}(l,i))'$$

$$= (0_{(l-i) \times i}, R_*^{(n)^{1/2}}(l-i)) R^{-1}(l) \mathbf{x}_j(l).$$

By an argument like that given in Lemma 4 of [17], we have, for all $1 \leq i < l \leq K_n$,

$$E \left\| \frac{1}{N} \sum_{j=K_n}^{n-\sqrt{n}-1} \mathbf{z}_j^{(n)}(l,i) e_{j+1} \right\|^{2q}$$

$$\leq C \left( \frac{l-i}{N} \right)^q \frac{1}{l-i} \sum_{r=1}^{l-i} \left( \frac{1}{N} \sum_{j=K_n}^{n-\sqrt{n}-1} E|z_{j,r}^{(n)}(l,i)|^{2q} \right),$$

and for all $1 \leq i < l \leq K_n$, all $1 \leq r \leq l-i$ and all $K_n \leq j \leq n - \sqrt{n} - 1$,
$$E|z_{j,r}^{(n)}(l,i)|^{2q} \leq C.$$
As a result, for all $1 \leq i < l \leq K_n$,

(5.60) $$E(|(\text{II})|^{2q}) \leq C\left(\frac{l-i}{N}\right)^q$$

follows.

Reasoning as for (5.59),

(5.61)
$$E(|(\text{III})|^{2q})$$
$$\leq CE\left\|\frac{1}{N}\sum_{j=K_n}^{n-\sqrt{n}-1}\mathbf{x}_j(i)e_{j+1}\right\|^{2q}_{(R^{-1}(i)-R_i^{-1}(l))'R_*^{(n)}(i)(R^{-1}(i)-R_i^{-1}(l))}$$

holds for all $1 \leq i < l \leq K_n$, where $R_i^{-1}(l)$ is the upper left $i \times i$ block of $R^{-1}(l)$, and $\mathbf{x}_j(i)$ and $R^{-1}(i)$, returning to their own original definitions, are an $i$-dimensional vector and an $i \times i$ matrix, respectively. If we write
$$R(l) = \begin{pmatrix} R(i) & R_{i,l-i}(l) \\ R_{l-i,i}(l) & R(l-i) \end{pmatrix},$$
then

(5.62) $$R_i^{-1}(l) = (R(i) - R_{i,l-i}(l)R^{-1}(l-i)R_{l-i,i}(l))^{-1}.$$

From (5.62),

(5.63) $$R^{-1}(i) - R_i^{-1}(l) = -R^{-1}(i)R_{i,l-i}(l)R^{-1}(l-i)R_{l-i,i}(l)R_i^{-1}(l).$$

On substituting (5.63) into the right-hand side of (5.61), an upper bound for the left-hand side of (5.61) is given by

(5.64) $$C\|M_n(l,i)\|^q E(\|\mathbf{u}_n(l,i)\|^{2q}),$$

where
$$M_n(l,i) = R_{l-i,i}(l)R^{-1}(i)R_*^{(n)}(i)R^{-1}(i)R_{i,l-i}(l)$$
and
$$\mathbf{u}_n(l,i) = R^{-1}(l-i)R_{l-i,i}(l)R_i^{-1}(l)\frac{1}{N}\sum_{j=K_n}^{n-\sqrt{n}-1}\mathbf{x}_j(i)e_{j+1}.$$

Observe that, for all $1 \leq i < l \leq K_n$,
$$\|M_n(l,i)\| \leq \|R^{-1/2}(i)R_*^{(n)}(i)R^{-1/2}(i)\|\|R_{l-i,i}(l)R^{-1}(i)R_{i,l-i}(l)\|$$
$$\leq C\|R_{l-i,i}(l)R^{-1}(i)R_{i,l-i}(l)\|$$
$$\leq C(\|R(l-i)\| + \|R(l-i) - R_{l-i,i}(l)R^{-1}(i)R_{i,l-i}(l)\|)$$
$$= C(\|R(l-i)\| + \|(R_{(l-i)^-}^{-1}(l))^{-1}\|)$$
$$\leq C,$$

where the second inequality follows from the boundedness of $\max_{1 \le i \le K_n} \|R_*^{(n)}(i)\|$, $R_{(l-i)^-}^{-1}(l)$ is the lower right $(l-i) \times (l-i)$ block of $R^{-1}(l)$, and the last inequality is ensured by the fact that, for all $1 \le i < l \le K_n$,

$$\|(R_{(l-i)^-}^{-1}(l))^{-1}\| \le \lambda_{\max}(R(K_n)) \le C,$$

where, for a symmetric matrix $A$, $\lambda_{\max}(A)$ denotes its maximum eigenvalue. As a result,

(5.65)  $$\max_{1 \le i < l \le k_n} \|M_n(l,i)\| \le C.$$

Moreover, by arguments similar to those used for showing (5.60) and (5.65), one has, for all $1 \le i < l \le K_n$,

$$E(\|\mathbf{u}_n(l,i)\|^{2q}) \le C \left(\frac{l-i}{N}\right)^q,$$

which, together with (5.61), (5.64) and (5.65), yields that, for all $1 \le i < l \le K_n$,

$$E(|(\mathrm{III})|^{2q}) \le C \left(\frac{l-i}{N}\right)^q.$$

This completes the proof of Lemma 8. $\square$

We are now ready to prove Theorem 2.

PROOF OF THEOREM 2.  By Hölder's inequality, one has, for $1 < r < \infty$,

(5.66)
$$E \left| \frac{\mathbf{f}_1(\hat{k}_n^S) - \mathbf{f}_1(k_n^*)}{L_n^{1/2}(\hat{k}_n^S)} \right|^{2q}$$
$$\le \sum_{k=1}^{K_n} \left( E \left| \frac{\mathbf{f}_1(k) - \mathbf{f}_1(k_n^*)}{L_n^{1/2}(k)} \right|^{2qr} \right)^{1/r} \{P(\hat{k}_n^S = k)\}^{(r-1)/r}.$$

Set $0 < \xi < \min\{1/2, \delta_1/2\}$. [Recall that $\xi$ is defined in (K.6) and $\delta_1$ is defined in (K.4).] Since (K.6) is assumed, there is a nonnegative number $\theta = \theta(\xi)$ such that (3.2) is fulfilled. By (3.2) and Lemma 8, the right-hand side of (5.66) is bounded by

(5.67)
$$C \left( \sum_{\substack{k=1 \\ k \notin A_{n,\theta}}}^{K_n} \left|\frac{k - k_n^*}{NL_n(k)}\right|^q + \sum_{\substack{k=1 \\ k \in A_{n,\theta}}}^{K_n} \left|\frac{k - k_n^*}{NL_n(k)}\right|^q \{P(\hat{k}_n^S = k)\}^{(r-1)/r} \right)$$
$$\le C \left( k_n^{*(\theta-1)q} k_n^{*\theta} + \sum_{\substack{k=1 \\ k \in A_{n,\theta}}}^{K_n} \left|\frac{k - k_n^*}{NL_n(k)}\right|^q \{P(\hat{k}_n^S = k)\}^{(r-1)/r} \right),$$

where $A_{n,\theta}$ is defined in (K.6). First observe that, for sufficiently large $q$, the first term on the right-hand side of (5.67) converges to 0. In addition, by analogy with (5.39) and the fact that for $a, b \geq 0$, $(a+b)^{(r-1)/r} \leq a^{(r-1)/r} + b^{(r-1)/r}$,

$$
\begin{aligned}
&\sum_{\substack{k=1 \\ k \in A_{n,\theta}}}^{K_n} \left| \frac{k - k_n^*}{NL_n(k)} \right|^q \{P(\hat{k}_n^S = k)\}^{(r-1)/r} \\
&\leq \sum_{i=1}^{5} \left\{ \sum_{\substack{k=1 \\ k \in A_{n,\theta}}}^{K_n} \left| \frac{k - k_n^*}{NL_n(k)} \right|^q \{P(|V_{in}(k)| > (1/5)V_n(k))\}^{(r-1)/r} \right\}.
\end{aligned}
\tag{5.68}
$$

By (3.2), (5.43) and (5.45), one has, for sufficiently large $n$ and $q$,

$$
\begin{aligned}
&\sum_{\substack{k=1 \\ k \in A_{n,\theta}}}^{K_n} \left| \frac{k - k_n^*}{NL_n(k)} \right|^q \{P(|V_{1n}(k)| > (1/5)V_n(k))\}^{(r-1)/r} \\
&\leq Ck_n^{*\xi q} \sum_{k=1}^{K_n} \left( \frac{k^q}{N^q} + N^{-q/2} \right) \\
&= o(1)
\end{aligned}
\tag{5.69}
$$

and

$$
\begin{aligned}
&\sum_{\substack{k=1 \\ k \in A_{n,\theta}}}^{K_n} \left| \frac{k - k_n^*}{NL_n(k)} \right|^q \{P(|V_{2n}(k)| > (1/5)V_n(k))\}^{(r-1)/r} \\
&\leq Ck_n^{*\xi q} \sum_{k=1}^{K_n} \left( \frac{k_n^{*q}}{N^q} + N^{-q/2} \right) \\
&= o(1).
\end{aligned}
\tag{5.70}
$$

According to (3.2) and (5.47), one has, for sufficiently large $n$ and $q$,

$$
\begin{aligned}
&\sum_{\substack{k=1 \\ k \in A_{n,\theta}}}^{K_n} \left| \frac{k - k_n^*}{NL_n(k)} \right|^q \{P(|V_{3n}(k)| > (1/5)V_n(k))\}^{(r-1)/r} \\
&\leq Ck_n^{*\xi q} \sum_{k=1}^{K_n} \left( \frac{k^q}{N^{q/2}} + \frac{k^{q/2}}{N^q L_n^q(k)} \right) \\
&\leq Ck_n^{*\xi q} \left( \frac{K_n^{q+1}}{N^{q/2}} + \sum_{k=1}^{k_n^*} k_n^{*-q/2} + \sum_{k=k_n^*+1}^{K_n} k^{-q/2} \right)
\end{aligned}
\tag{5.71}
$$

$$= o(1),$$

where the last equality is ensured by $0 < \xi < \min\{1/2, \delta_1/2\}$. Similarly, by (3.2) and (5.48),

$$(5.72) \quad \sum_{\substack{k=1 \\ k \in A_{n,\theta}}}^{K_n} \left|\frac{k - k_n^*}{NL_n(k)}\right|^q \{P(|V_{4n}(k)| > (1/5)V_n(k))\}^{(r-1)/r} = o(1),$$

provided $q$ is sufficiently large. Finally, by (3.2) and (5.49), one has, for sufficiently large $n$ and $q$,

$$\sum_{\substack{k=1 \\ k \in A_{n,\theta}}}^{K_n} \left|\frac{k - k_n^*}{NL_n(k)}\right|^q \{P(|V_{5n}(k)| > (1/5)V_n(k))\}^{(r-1)/r}$$

$$(5.73) \quad \leq C k_n^{*\xi q} \sum_{k=1}^{K_n} \frac{\|\mathbf{a}(k) - \mathbf{a}(k_n^*)\|_R^q}{N^{q/2} L_n^q(k)}$$

$$\leq C k_n^{*\xi q} \left(\sum_{k=1}^{k_n^*} k_n^{*-q/2} + \sum_{k=k_n^*+1}^{K_n} k^{-q/2}\right)$$

$$= o(1).$$

In view of (5.69)–(5.73), for sufficiently large $q$, the left-hand side of (5.68) converges to 0. Consequently, the left-hand side of (5.67) also converges to 0 for sufficiently large $q$. This result, Jensen's inequality and (5.66) yield

$$\lim_{n \to \infty} E\left|\frac{\mathbf{f}_1(\hat{k}_n^S) - \mathbf{f}_1(k_n^*)}{L_n^{1/2}(\hat{k}_n^S)}\right|^{2q} = 0,$$

for any $q > 0$. Moreover, by Theorem 3 and Lemma 7,

$$(5.74) \quad \lim_{n \to \infty} E\left(\frac{|\mathbf{f}(\hat{k}_n^S) - \mathbf{f}(k_n^*)|^{2q}}{L_n^q(k_n^*)}\right) = 0.$$

On the other hand, by Wei ([29], Lemma 2) it can be shown that

$$E|\mathcal{S}(k) - \mathcal{S}(k_n^*)|^{2q} \leq C \|\mathbf{a}(k) - \mathbf{a}(k_n^*)\|_R^{2q}.$$

By the definitions of $L_n(k)$ and $k_n^*$, it is easy to see that

$$\|\mathbf{a}(k) - \mathbf{a}(k_n^*)\|_R^{2q} \leq \left(L_n(k) - L_n(k_n^*) + \left|\frac{k - k_n^*}{N}\right|\sigma^2\right)^q.$$

These facts and an argument similar to that used for verifying (5.74) yield

$$(5.75) \quad \lim_{n \to \infty} E\left(\frac{|\mathcal{S}(\hat{k}_n^S) - \mathcal{S}(k_n^*)|^{2q}}{L_n^q(k_n^*)}\right) = 0.$$

As a result, (3.8) with $\hat{k}_n = \hat{k}_n^S$ follows from (5.74), (5.75), Proposition 2 and the Cauchy–Schwarz inequality. Moreover, by arguments similar to those used for showing (5.74), (5.75) and Corollary 2, (3.8) also holds with $\hat{k}_n = \hat{k}_n^A, \hat{k}_n^F, \hat{k}_n^{S_p}$ or $\hat{k}_n^C$. $\square$

PROOF OF COROLLARY 1. The proof of Corollary 1 is similar to that of Theorem 2. The details are omitted. $\square$

**6. Discussion and concluding remarks.** (1) Main contributions. Due to its success in practical applications, AIC has received considerable attention among researchers (and practitioners) from various disciplines in the past three decades; see [10]. However, the statistical properties of AIC are still not clear when it is used for forecasting the future of an observed time series. In the present article we have attempted to resolve this problem. Armed with some new technical tools, Theorem 2 and (4.4) successfully show how well AIC (and its variants) can work for same-realization predictions. The simulation results given in Table 1 also show that the finite-sample performance of AIC is satisfactory in many practical situations. Corollary 1 and (3.12)–(3.14) explore the prediction efficiencies of some other AIC-like criteria having different penalties for the number of regressors in the model. By the same argument used for proving Theorem 2 and the ideas of Shibata ([27], Section 5), Bhansali [5] and Ing and Wei ([17], Section 4), extending Theorem 2 and Corollary 1 to the multistep prediction case is straightforward.

On the other hand, Table 2 indicates that for same-realization predictions it seems very difficult for AIC to be strongly asymptotically efficient. This phenomenon not only points out another dissimilarity between same- and independent-realization predictions [since AIC is strongly asymptotically efficient for independent-realization predictions; see (4.10)], it also inspires a new direction for time series model selection, that is, selecting models (or orders) through the second-order conditional MSPE, namely, $E[(x_{n+1} - \hat{x}_{n+1}(k))^2 | x_1, \ldots, x_n] - \sigma^2$. As suggested by Table 2, a model having the minimal second-order conditional MSPE can be asymptotically much more efficient than the model selected by AIC in the sense that $r_n^* \gg 1$ for all large $n$. [Recall that $r_n^*$ is defined after (4.8).] Unfortunately, since $E[(x_{n+1} - \hat{x}_{n+1}(k))^2 | x_1, \ldots, x_n] - \sigma^2$ is unobservable, it cannot be used as a selection criterion in practice. However, we conjecture that a model selection criterion based on a reliable estimator of $E[(x_{n+1} - \hat{x}_{n+1}(k))^2 | x_1, \ldots, x_n] - \sigma^2$ should also outperform AIC for same-realization predictions. For a related discussion on estimating $E[(x_{n+1} - \hat{x}_{n+1}(k))^2 | x_1, \ldots, x_n] - \sigma^2$ in finite-order AR models, see [17].

(2) Moment restrictions. For independent-realization predictions, we provide a set of sufficient conditions in Remark 6 which guarantees that AIC achieves Shibata's asymptotic efficiency in non-Gaussian settings. Since Remark 6 only assumes that the error distributions have uniformly bounded $(4+\delta_0)$th moments, where $\delta_0$ is any (small) positive number, it notably improves the best previous result given by Lee and Karagrigoriou ([21], Theorem 3.1), which required the existence of the eighth moment. However, to ensure that AIC is asymptotically efficient in the senses of (4.3) and (4.4), (K.3) is needed; see Theorems 1 and 2. This is because (2.4) and (2.5) are required to hold for sufficiently large $q$ in the proofs of Theorems 1 and 2; and Proposition 1, which is the first result giving sufficient conditions [including (K.3)] such that (2.4) and (2.5) are fulfilled for any $q > 0$, is used to meet this requirement. Although (K.3) is rather stronger than is necessary, it is convenient. Note that if (K.4) is replaced by a stronger assumption than $K_n^{6+\delta_0} = O(n)$ for some $\delta_0 > 0$, then the first part of Theorem 2 of [17] holds instead of Proposition 1. By applying this result to our analysis, (K.3) in Proposition 2 and Theorem 1 can be weakened to $\sup_{-\infty < t < \infty} E|e_t|^{20} < \infty$ and $\sup_{-\infty < t < \infty} E|e_t|^{36} < \infty$, respectively. But for Theorem 2, (K.3) needs to be replaced by a more complicated moment restriction which may depend on the value of $\theta$ [which is defined in (K.6)]. Consequently, while (K.3) can be slightly relaxed at the price of reducing the number of candidate models, the disadvantages seem to outweigh the merits. To resolve this dilemma, it is necessary to verify (2.4) and (2.5) under milder moment conditions. However, this topic is beyond the scope of the present article.

(3) *Extensions to the regression model.* As mentioned in Section 1, the main difficulty of analyzing AIC's (and its variants') same-realization MSPE in time series models lies in the fact that future observations, estimated parameters and selection criteria are all stochastically dependent and, hence, Shibata's approach is no longer applicable. For the regression model, however, the (commonly used) assumption of independent observations yields that the future observations are independent of the estimated parameters and selection criteria even in the same-realization case, which substantially simplifies the task of analyzing the model selection criterion's (same-realization) MSPE. This is exactly the same situation encountered in independent-realization predictions for time series models (see Section 1). In fact, under a Gaussian regression model with infinitely many parameters, Breiman and Freedman [7] showed that the $S_p$ is asymptotically efficient for "same-realization" predictions from a conditional MSPE point of view. [Note that their asymptotic efficiency is the same as the one discussed in (4.11).] It also can be shown that their result still holds with $S_p$ replaced by AIC, FPE, $C_p$ or $S_n(k)$. In addition, we conjecture that their result can be extended to unconditional versions without the Gaussian assumption, provided suitable smoothness conditions [such as (K.2)] on the distributions of the (random) regressors and the white noise are imposed.

(4) The possibility of extensions to the multivariate case. Since for multivariate time series, AR models are the most used models by far, order selection for vector AR models has attracted growing interest among researchers from various disciplines in recent years. For example, Findley and Wei [12] recently presented the first mathematically complete derivation of the multivariate AIC for comparing vector AR models fit to stationary series in independent-realization settings. When the candidate vector AR models are misspecified, Schorfheide [23] also considered order selection problems for the purpose of independent-realization predictions. However, since these results focus on independent-realization cases, it would be of interest to extend our same-realization results to the multivariate case. To achieve this goal, extending Proposition 1 to stationary multivariate time series is necessary. Taking the approaches used to verify Theorem 4.1 of [12] [which gives a multivariate version of (2.4) but with $K_n$ fixed with $n$] and Theorem 2 of [17], a multivariate extension of Proposition 1 can be easily obtained. In addition, generalizations of the moment bounds of Section 5 to the vector case are also required for establishing the desired results. However, since these generalizations are not straightforward, further investigation along this direction is needed.

## APPENDIX

PROOF OF (3.4). We first show that, for sufficiently large $n$,

(A.1) $$\frac{1}{\beta} \log N - C \log_2 N \leq k_n^* \leq \frac{1}{\beta} \log N + C \log_2 N,$$

where $C$ is some positive number and $\beta$ is defined in (3.3). Let $k = k_n^* - d$, where $d$ is some positive integer. Then

(A.2) $$L_n(k) - L_n(k_n^*) = \|\mathbf{a} - \mathbf{a}(k)\|_R^2 - \|\mathbf{a} - \mathbf{a}(k_n^*)\|_R^2 - \frac{d}{N}\sigma^2 \geq 0.$$

According to (A.2) and (3.3), one has, for some $C > 0$,

$$\frac{d}{N}\sigma^2 \leq C k_n^{*\theta_1} e^{-\beta(k_n^* - d)},$$

where $\theta_1$ is defined in (3.3). Taking the (natural) logarithm of both sides, we get

(A.3) $$k_n^* \leq \frac{1}{\beta} \log N + C \log_2 N,$$

for some $C > 0$. In view of (A.3) and (K.4), we have $k_n^* + k_n^{*\eta} \leq K_n$ for any $0 < \eta < 1$ and for sufficiently large $n$. Now let $k = k_n^* + k_n^{*\eta}$ for some $0 < \eta < 1$.

(Here we assume without loss of generality that $k_n^{*\eta}$ is an integer.) In this case

$$(\text{A.4}) \quad L_n(k) - L_n(k_n^*) = \frac{k_n^{*\eta}}{N}\sigma^2 - (\|\mathbf{a} - \mathbf{a}(k_n^*)\|_R^2 - \|\mathbf{a} - \mathbf{a}(k)\|_R^2) \geq 0.$$

By (3.3), (A.4) and the fact that $k_n^* \to \infty$, one has, for sufficiently large $n$ and some $C > 0$,

$$\frac{k_n^{*\eta}}{N}\sigma^2 \geq C k_n^{*-\theta_1} e^{-\beta k_n^*}$$

and, hence,

$$(\text{A.5}) \qquad k_n^* \geq \frac{1}{\beta}\log N - C \log_2 N.$$

Consequently, (A.1) follows from (A.3) and (A.5).

To show (3.4), first assume that $k < k_n^* - k_n^{*\eta}$ for some $0 < \eta < 1$. By (3.3) one has, for sufficiently large $n$ and some $C > 0$,

$$(\text{A.6}) \quad \begin{aligned}\frac{L_n(k) - L_n(k_n^*)}{(k_n^* - k)/N} &= \frac{\|\mathbf{a} - \mathbf{a}(k)\|_R^2 - \|\mathbf{a} - \mathbf{a}(k_n^*)\|_R^2 - ((k_n^* - k)/N)\sigma^2}{(k_n^* - k)/N} \\ &\geq \frac{C k_n^{*-\theta_1} e^{-\beta(k_n^* - k_n^{*\eta})}}{k_n^*/N} - \sigma^2.\end{aligned}$$

In view of (A.1), the right-hand side of (A.6) diverges to infinity and, hence, (3.4) holds for $k < k_n^* - k_n^{*\eta}$.

For $k > k_n^* + k_n^{*\eta}$, $0 < \eta < 1$, one has

$$(\text{A.7}) \quad \begin{aligned}\frac{L_n(k) - L_n(k_n^*)}{(k - k_n^*)/N} &\geq \frac{L_n(k) - L_n(k_n^* + (1/2)k_n^{*\eta})}{(k - k_n^*)/N} \\ &\geq \frac{\sigma^2}{2} - \frac{\|\mathbf{a} - \mathbf{a}(k_n^* + (1/2)k_n^{*\eta})\|_R^2}{k_n^{*\eta}/N}.\end{aligned}$$

By (3.3) and (A.1), the second term on the right-hand side of (A.7) converges to 0. Therefore, (3.4) holds for $k > k_n^* + k_n^{*\eta}$. □

PROOF OF (3.7). Let $k = k_n^* - d$, where $d$ is some positive integer. By (3.5) and (A.2),

$$(\text{A.8}) \quad \frac{d\sigma^2}{N} \leq (C_4 + M_1(k_n^* - d)^{-\xi_1})(k_n^* - d)^{-\beta} - (C_4 - M_1 k_n^{*-\xi_1})k_n^{*-\beta}.$$

By Taylor's theorem, (A.8) can be further expressed as

$$\frac{d\sigma^2}{N} \leq C_4 \beta d k_n^{*-1-\beta} + O(k_n^{*-2-\beta}).$$

Therefore, for sufficiently large $n$ and some $C > 0$,

$$\left(\frac{\sigma^2}{NC_4\beta}\right)^{-1/(\beta+1)} \geq k_n^* - C.$$

Since $\beta > 1 + \delta_1$, we can choose a positive integer, $d$, such that $k_n^* + d \leq K_n$ for all large $n$. By letting $k = k_n^* + d$ and by an argument analogous to that used in (A.8), we have, for sufficiently large $n$ and some $C > 0$,

$$\left(\frac{\sigma^2}{NC_4\beta}\right)^{-1/(\beta+1)} \leq k_n^* + C.$$

As a result,

(A.9) $$k_n^* = \left(\frac{\sigma^2}{NC_4\beta}\right)^{-1/(\beta+1)} + O(1).$$

Armed with (A.9), the proof of (3.7) is divided into four cases.

*Case* 1. $1 \leq k \leq \theta_2 k_n^*$, where $0 < \theta_2 < 1$ is chosen to satisfy $\theta_2^{-\beta} - 1 > \beta$. By (3.5) and (A.9), one has

$$L_n(k) - L_n(k_n^*) = \|\mathbf{a} - \mathbf{a}(k)\|_R^2 - \|\mathbf{a} - \mathbf{a}(k_n^*)\|_R^2 - \frac{k_n^* - k}{N}\sigma^2$$

$$\geq \|\mathbf{a} - \mathbf{a}(\theta_2 k_n^*)\|_R^2 - \|\mathbf{a} - \mathbf{a}(k_n^*)\|_R^2 - \frac{k_n^*}{N}\sigma^2$$

$$= \frac{k_n^*}{N\beta}\sigma^2(\theta_2^{-\beta} - 1 - \beta + o(1)).$$

Therefore,

(A.10) $$\frac{L_n(k) - L_n(k_n^*)}{(k_n^* - k)/N} > C$$

holds for sufficiently large $n$, $1 \leq k \leq \theta_2 k_n^*$, and some $C > 0$.

*Case* 2. $\theta_2 k_n^* < k < k_n^*$, where $\theta_2$ is defined as in case 1.
By Taylor's theorem, (3.5) and the assumption that $\xi_1 \geq 2$,

$$\frac{L_n(k) - L_n(k_n^*)}{(k_n^* - k)/N}$$

(A.11) $$= \left(\frac{d}{N}\right)^{-1}\left\{\left[k_n^{*-\beta}(C_4 + O(k_n^{*-\xi_1}))\right.\right.$$
$$\left.\times \left(1 + \frac{\beta d}{k_n^*} + \frac{\beta(\beta+1)d^2}{2k_n^{*2}}\right.\right.$$

$$+\frac{\beta(\beta+1)(\beta+2)d^3}{6k_n^{*3}}(1-\alpha_n)^{-\beta-3}\bigg)$$

$$-k_n^{*-\beta}(C_4+O(k_n^{*-\xi_1}))\bigg]-\frac{d}{N}\sigma^2\bigg\}$$

$$=\left(\frac{d}{N}\right)^{-1}\left[\frac{C_4\beta d}{k_n^{*\beta+1}}+\frac{C_4 d^2\beta\gamma_{\beta,\theta_2}}{k_n^{*2+\beta}}-\frac{d}{N}\sigma^2+O(k_n^{*-2-\beta})\right],$$

where $d=k_n^*-k$, $0<\alpha_n<dk_n^{*-1}\leq 1-\theta_2$ and

$$\gamma_{\beta,\theta_2}=\frac{\beta+1}{2}+\frac{(\beta+1)(\beta+2)d}{6k_n^*(1-\alpha_n)^{\beta+3}}.$$

In view of (A.9) and (A.11), we have

(A.12)
$$\frac{L_n(k)-L_n(k_n^*)}{(k_n^*-k)/N}=\frac{\sigma^2}{k_n^*}(d\gamma_{\beta,\theta_2}+O(1))$$
$$\geq C\frac{d}{k_n^*}$$

for sufficiently large $n$ and some $C>0$, provided $d>C_6$ with some $C_6>0$.

*Case* 3. $(1+\theta_2)k_n^*<k<K_n, 0<\theta_2<1$.
By (3.5) and (A.9),

(A.13)
$$L_n(k)-L_n(k_n^*)=\frac{d}{N}\sigma^2-C_4 k_n^{*-\beta}(1-(1+(d/k_n^*))^{-\beta})$$
$$+O(k_n^{*-\xi_1-\beta})$$
$$=\frac{d}{N}\sigma^2-\frac{k_n^*}{\beta N}\sigma^2(1-(1+(d/k_n^*))^{-\beta})(1+o(1))$$
$$+O(k_n^{*-\xi_1-\beta}),$$

where $d=k-k_n^*$. By (A.13) one has, for sufficiently large $n$,

(A.14)
$$\frac{L_n(k)-L_n(k_n^*)}{(k-k_n^*)/N}=\sigma^2\left(1-\frac{1}{\beta(d/k_n^*)}(1-(1+(d/k_n^*))^{-\beta})(1+o(1))\right)$$
$$+O\left(\frac{N}{d}k_n^{*-\xi_1-\beta}\right)$$
$$>C>0,$$

where the last inequality follows from

$$O\left(\frac{N}{d}k_n^{*-\xi_1-\beta}\right)=o(1)$$

and

$$g(x) = 1 - \frac{1}{\beta x}(1 - (1+x)^{-\beta}) > 0,$$

for $x > 0$.

*Case* 4. $k_n^* < k \leq (1+\theta_2)k_n^*$ with $0 < \theta_2 < 3/(\beta+2)$.

By Taylor's theorem, (3.5), (A.9) and the assumption that $\xi_1 \geq 2$, one has, for sufficiently large $n$ and some $C > 0$,

$$\frac{L_n(k) - L_n(k_n^*)}{(k-k_n^*)/N}$$
$$= \left(\frac{d}{N}\right)^{-1}$$
$$\times \left\{\frac{d}{N}\sigma^2 - \frac{\sigma^2 k_n^*}{N\beta}(1 + O(N^{-1/(\beta+1)}))\right.$$
$$\left. \times \left(\frac{\beta d}{k_n^*} - \frac{\beta(\beta+1)d^2}{2k_n^{*2}}\right.\right.$$
$$\left.\left. + \frac{\beta(\beta+1)(\beta+2)d^3}{6k_n^{*3}}(1+\alpha_n)^{-\beta-3}\right) + O(k_n^{*-2-\beta})\right\}$$

(A.15)

$$= \sigma^2 \left[\frac{(\beta+1)d}{2k_n^*} - \frac{(\beta+1)(\beta+2)d^2}{6k_n^{*2}}(1+\alpha_n)^{-3-\beta} + O(k_n^{*-1})\right]$$
$$> C\frac{d}{k_n^*},$$

provided $d = k - k_n^* > C_6$ for some $C_6 > 0$. Here $\alpha_n$ is some positive number which satisfies $0 < \alpha_n < d/k_n^*$, and the last inequality follows from the condition on $\theta_2$.

Consequently, the desired property (3.7) is ensured by (A.10), (A.12), (A.14) and (A.15). $\square$

**Acknowledgments.** We are deeply grateful to Dr. David Findley for helpful discussions and suggestions on this work. We also would like to thank a Co-Editor, an Associate Editor and three anonymous referees for their valuable and constructive comments which helped to improve the paper.

INSTITUTE OF STATISTICAL SCIENCE
ACADEMIA SINICA
DEPARTMENT OF ECONOMICS
NATIONAL TAIWAN UNIVERSITY
TAIPEI 11529
TAIWAN
REPUBLIC OF CHINA
E-MAIL: cking@stat.sinica.edu.tw

INSTITUTE OF STATISTICAL SCIENCE
ACADEMIA SINICA
DEPARTMENT OF MATHEMATICS
NATIONAL TAIWAN UNIVERSITY
TAIPEI 11529
TAIWAN
REPUBLIC OF CHINA